\documentclass{article}%
\usepackage{amsmath}
\usepackage{amsfonts}
\usepackage{eurosym}
\usepackage{amssymb}
\usepackage{graphicx}%
\providecommand{\U}[1]{\protect\rule{.1in}{.1in}}
\newtheorem{theorem}{Theorem}

\newtheorem{corollary}[theorem]{Corollary}

\newtheorem{definition}[theorem]{Definition}

\newtheorem{lemma}[theorem]{Lemma}

\newenvironment{proof}[1][Proof]{\noindent\textbf{#1.} }{\ \rule{0.5em}{0.5em}}
\begin{document}

\title{On the Kneser property of the three-dimensional Navier-Stokes equations with damping}
\author{D. Pardo, J. Valero, A. Gim\'{e}nez\\{\small Universidad Miguel Hern\'{a}ndez de Elche,. Centro de
Investigaci\'{o}n Operativa} \\{\small Avda. Universidad s/n, 03202-Elche (Alicante), Spain}\\{\small E.mails: daniel.pardo@alu.umh.es, a.gimenez@umh.es, jvalero@umh.es}}
\date{}
\maketitle

\begin{abstract}
In this paper, we study the connectedness and compactness of the attainability
set of weak solutions to the three-dimensional Navier--Stokes equations with
damping. Depending on the value of the parameter $\beta$, which controls the
damping term, we establish these results with respect to either the weak or
the strong topology of the phase space. In the latter case, we also prove that
the global attractor is connected.

Additionally, we establish results concerning the regularity of the global
attractor and provide a new proof of its existence for strong solutions.

\end{abstract}



\bigskip

\section{Introduction}

In this paper, we aim to study the Kneser property---that is, the
connectedness and compactness---of the attainability set of solutions to the
three-dimensional Navier--Stokes equations with damping, given by%
\begin{equation}
\left\{
\begin{array}
[c]{c}%
u_{t}-\mu\Delta u+(u\cdot\nabla)u+\nabla p+\alpha|u|^{\beta-1}u=f,\ (x,t)\in
\Omega\times(\tau,T),\\
\operatorname{div}u=0,\ \\
u|_{\partial\Omega}=0,\\
u(\tau)=u_{\tau},
\end{array}
\right.  \label{eq:NSDE-1-intro}%
\end{equation}
where $\alpha,\mu>0$ and $\beta\geq1$ are constants, $\Omega$ is a bounded
open subset of $\mathbb{R}^{3}$ with smooth boundary $\partial\Omega$. The
functions $u(x,t)=(u_{1}(x,t),$ $u_{2}(x,t),$ $u_{3}(x,t))$, $p(x,t)$ stand
for the velocity field and the pressure, respectively. This system describes
the evolution of an incompressible fluid with resistance to the motion. It is
appropriate to model various processes, such as thermal dispersion in a porous
media \cite{HsuCheng} or the flow of cerebrospinal fluid inside the porous
brain tissues \cite{bib12}. Additionally, compressible Euler equations with
damping describe the flow of a compressible gas through a porous medium
\cite{bib08}, whereas Saint-Venant equations are used in oceanography to
describe the flow of viscous shallow water with friction \cite{bib01}. These
equations are known as Brinkman-Forchheimer equations or tamed Navier-Stokes equations.

There is extensive literature on this problem. In the work \cite[Theorem
1.2]{bib02}, the existence of global weak solutions to problem
(\ref{eq:NSDE-1-intro}) with initial data in the space of square-integrable,
divergence-free functions was established for $\beta\geq1$ and $\Omega
=\mathbb{R}^{3}$. This result was extended later in \cite[Theorem 2.1]{bib15}
to bounded domains with Dirichlet boundary conditions. The uniqueness of weak
solutions was established for bounded domains in \cite{bib32} for $\beta\geq
4$, and subsequently extended to the cases $\beta>3$ and $\beta=3$,
$4{\alpha\mu}\geq1$ in \cite{bib26}. For periodic boundary conditions the
existence and uniqueness of weak solutions was established for $\beta>3$ and
$\beta=3,\ 2{\alpha\mu}\geq1$ \cite{KinraMoan}$.$

Regarding global strong solutions with $\Omega$ bounded and Dirichlet boundary
conditions, their existence was proved in \cite{bib27} for $\beta\in\left(
3,5\right)  $ and $\beta=3,\ 4{\alpha\mu}>1$ and for $\beta=5$ in
\cite{LiKimKimO}. In the previous papers \cite{bib15}, \cite{bib17} the
existence of strong solutions was stated to be true for $\beta>3$, but it is
pointed out in \cite{bib27}, \cite{LiKimKimO} that the proof in those works is
unclear as the function $-\Delta u$ is used in quality of test function, which
seems to be incorrect for $\Omega\not =\mathbb{R}^{3}$ with Dirichlet boundary
conditions. For periodic boundary conditions the existence of global strong
solutions was established for $\beta>3$ and $\beta=3,\ 4{\alpha\mu\geq}1$ in
\cite{HajduRob} (see also \cite{bib34}, \cite{GautamMohan}). We observe that
these results are stronger because on domains with periodic conditions it is
right to use $-\Delta u$ as test function.

Other results on existence, uniqueness and regularity of local and global
strong solutions can be found in \cite{bib02}, \cite{LiKimKimO}, \cite{Tinh},
\cite{bib22}, \cite{bib24}, \cite{bib23}, \cite{bib24}, \cite{bib25}.

The asymptotic behavior of weak solutions to problem (\ref{eq:NSDE-1-intro})
was studied in \cite{bib32}. More precisely, for $\beta\geq3$ a multivalued
semiflow consisting of weak solutions was defined and the existence of a
global compact attractor was proved. It was stated in \cite{bib32} that the
attractor is more regular for $\beta\geq4,$ but the proof given in that paper
is unclear because, as in \cite{bib15}, \cite{bib17}, the function $-\Delta u$
is used in quality of test function. A different proof for either $\beta
\in\left(  3,5\right)  $ or $\beta=3$, $4{\alpha\mu}>1$ is given in the
Appendix. Thus, this question remains open for $\beta\geq5.$ For periodic
boundary conditions and either $\beta>3$ or $\beta=3,\ 2{\alpha\mu}\geq1$, in
\cite{KinraMoan} a semigroup of operator is defined, proving the existence of
the global attractor.

In \cite{Tinh} the existence of the global attractor for weak solutions was
proved in the situation where $-\Delta u$ is replaced by the fractional
laplacian and the damping term is exponential. Additionally, the convergence
to $0$ of weak solutions as time goes to infinity in the particular case where
$f=0$ and $\Omega=\mathbb{R}^{3}$ has been studied in several papers (see
\cite{Peng} and its references).

The asymptotic behavior of strong solution to problem (\ref{eq:NSDE-1-intro})
was studied in \cite{bib15}-\cite{bib17} in the autonomous and non-autonomous
settings. More precisely, the existence of the global attractor was stated for
$\beta\in(3,5]$. Again, the proof of this result is unclear because the
function $-\Delta u$ was used in quality of test function. We give a new proof
of the existence of the strong attractor for $\beta\in\left(  3,5\right)  $
and $\beta=3$, $4{\alpha\mu}>1$. For $\beta\geq5$ this question remains open.
For periodic boundary conditions and either $\beta\in(3,5]$ or $\beta
=3,\ 2{\alpha\mu}\geq1$ this result is proved in \cite{KinraMoan}, and for
$\beta=3,\ 4{\alpha\mu}\geq1$ in \cite{HajduRob} but restricting the semigroup
to a suitable ball.

When uniqueness of the Cauchy problem fails, an interesting question is to
study the connectedness and compactness of the attainability set of solutions
at every future moment of time. This is known as the Kneser property in the
literature. In \cite{bib30}, such result was established for\ weak solutions
to the three-dimensional Navier-Stokes system with respect to the weak
topology of the phase space. In this paper, we prove the Kneser property for
weak solutions to problem (\ref{eq:NSDE-1-intro}) with respect to the weak
topology if $1\leq\beta<3$ and to the strong topology if $\beta=3$,
$4{\alpha\mu}<1$. To this end, we introduce as an auxiliary tool the so-called
globally modified Navier-Stokes equations with damping and prove several
properties of their solutions. In addition, using the Kneser property we prove
that the global attractor for weak solutions is connected with respect to the
strong tolopology of the phase space for $\beta=3,\ 4\alpha\mu<1.$

This paper is organized as follows.

In Section 2, we recall several definitions and results concerning solutions
to problem (\ref{eq:NSDE-1-intro}). In Section 3, we introduce the Globally
Modified Navier-Stokes Equations with Damping and prove the existence,
uniqueness and continuity of the weak solutions. In Section 4, we establish
the Kneser property of weak solutions with respect to the weak topology of the
phase space for $\beta\in\lbrack1,3]$. In Section 5, we demonstrate the Kneser
property of weak solutions with respect to the strong topology of the phase
space for $\beta=3$, $4{\alpha\mu}<1$. In Section 6, we show that the global
attractor of weak solutions is connected with respect to the strong topology
of the phase space for $\beta=3$, $4{\alpha\mu}<1$ and prove the existence of
the strong attractor for $\beta\in\left(  3,5\right)  $ and $\beta=3$,
$4{\alpha\mu}>1$. In Section 7, we prove the regularity of the weak attractor
for $\beta\in\left(  3,5\right)  $ and $\beta=3$, $4{\alpha\mu>}1.$

\section{Preliminaries}

This section recalls some results concerning the existence and properties of
weak solutions to problem (\ref{eq:NSDE-1-intro}).

Here and further $|\cdot|_{d}$ denotes in general the norm in $\mathbb{R}^{d}$
for any $d\geq1$. We define the usual function spaces
\begin{align*}
\mathcal{V}  &  =\{u\in(C_{0}^{\infty}(\Omega))^{3}:\operatorname{div}u=0\},\\
H  &  =cl_{(L^{2}(\Omega))^{3}}\mathcal{V},\\
V  &  =cl_{(H_{0}^{1}(\Omega))^{3}}\mathcal{V},
\end{align*}
where $cl_{X}$ denotes the closure in the space $X$. It is well known that
$H$, $V$ are separable Hilbert spaces and identifying $H$ and its dual we have
$V\subset H\subset V^{\prime}$ with dense and continuous injections. We denote
by $\left(  \text{\textperiodcentered,\textperiodcentered}\right)  $,
$|\cdot|$ and $\left(  \left(  \text{\textperiodcentered,\textperiodcentered
}\right)  \right)  $, $||\cdot||$, the inner product and norm in $H$ and $V$ ,
respectively, and by $\langle\cdot,\cdot\rangle$ the duality between
$V^{\prime}$ and $V$, while $|\cdot|_{q},\left(  \text{\textperiodcentered
,\textperiodcentered}\right)  $ denote the norm and duality in $(L^{q}%
(\Omega))^{3}$ for $q>2$. Let $H_{w}$ be the space $H$ endowed with the weak
topology. As usual,we define the continuous trilinear form $b:V\times V\times
V\rightarrow\mathbb{R}$ by
\[
b(u,v,w)=\sum_{i,j=1}^{3}{u_{i}\frac{\partial v_{j}}{\partial x_{i}}w_{j}}.
\]
It is well-known that $b(u,v,v)=0$, if $u\in V$, $v\in(H_{0}^{1}(\Omega))^{3}%
$. For $u,v\in V$ we denote by $B(u,v)$ the element of $V^{\prime}$ such that
$\langle B(u,v),w\rangle=b(u,v,w)$, for all $w\in V$. Let $P$ be the
orthogonal projector from $(L^{2}(\Omega))^{3}$ onto $H$ and {$Au=-P\Delta u$}
be the Stokes operator, defined by $\langle Au,v\rangle=((u,v))$ for $u,v\in
V$. Since the boundary $\partial\Omega$ is smooth, $D(A)=(H^{2}(\Omega
))^{3}\cap V$ and $|Au|$ defines a norm in $D(A)$ which is equivalent to the
norm in $(H^{2}(\Omega))^{3}$. Also, we define $G(u)$ by $G(u)=\alpha
P|u|^{\beta-1}u$.

\begin{definition}
\label{WeakSol}For $u_{\tau}\in H,\ f\in L^{2}(\tau,T;V^{\prime})$, $\beta
\geq1$ the function
\[
u\in L^{\infty}(\tau,T;H)\cap L^{2}(\tau,T;V)\cap L^{\beta+1}(\tau
,T;(L^{\beta+1}(\Omega))^{3})
\]
is said to be a weak solution to problem (\ref{eq:NSDE-1-intro}) on $[\tau,T]$
if $u(\tau)=u_{\tau}$ and
\begin{equation}
\frac{d}{dt}(u,v)+\mu((u,v))+b(u,u,v)+\alpha(|u|^{\beta-1}u,v)=\langle
f,v\rangle, \label{eq:NSDE-3-weaksol}%
\end{equation}
for any $v\in V\cap(L^{\beta+1}(\Omega))^{3}$, in the sense of scalar distributions.
\end{definition}

For each $T>\tau$ if $u$ satisfies equation \eqref{eq:NSDE-3-weaksol}, then%
\[
u\in C([0,T],H_{w}),\frac{du}{dt}\in L^{\frac{4}{3}}(\tau,T;V^{\prime
})+L^{\frac{\beta+1}{\beta}}(\tau,T;(L^{\frac{\beta+1}{\beta}}(\Omega
))^{3})\text{ for all }T>\tau.
\]

In particular, the initial condition $u(\tau)=u_{\tau}$ makes sense for any
$u_{\tau}\in H$ and the equation \eqref{eq:NSDE-3-weaksol} is equivalent to
\[
(u(t),w)+\mu\int_{\tau}^{t}(u(s),w)ds+\int_{\tau}^{t}b(u(s),u(s),w)ds+\alpha
\int_{\tau}^{t}(|u(s)|^{\beta-1}u,w)=(u_{\tau},w)+\int_{\tau}^{t}\langle
f(s),w\rangle ds,
\]
for all $t>\tau$ and $w\in V\cap L^{\beta+1}(\Omega)$. It is known
\textbf{\cite{bib16} } that if $f\in L^{2}(\tau,T;V^{\prime})$ and $u_{\tau
}\in H,$ then problem (\ref{eq:NSDE-1-intro}) has at least one weak solution.
Moreover, as this solution is obtained via Galerkin approximations, we can
prove that it satisfies the energy inequality
\begin{equation}
V_{\tau}(u(t))\leq V_{\tau}(u(s)),\ \forall s\leq t\leq T,\ \text{for a.e.
}s>\tau\text{ and }s=\tau, \label{eq:NSDE-4-Eineq}%
\end{equation}
where
\[
V_{\tau}(u(t)):=\frac{1}{2}\left\vert u(t)\right\vert ^{2}+\mu\int_{\tau}%
^{t}||u(r)||^{2}+\alpha\int_{\tau}^{t}|u_{n}(r)|_{\beta+1}^{\beta+1}%
dr-\int_{\tau}^{t}\langle f(r),u(r)\rangle dr.
\]
Indeed, if $u_{n}$ are the Galerkin approximations, then%
\[
V_{\tau}(u_{n}(t))=V_{\tau}(u_{n}(s))\text{ for all }t\geq s\geq\tau,
\]
and passing to the limit in a standard way we get (\ref{eq:NSDE-4-Eineq}).

By the same argument as in the proof from \cite[Lemma 3.8]{bib32} we obtain
that the translation and concatenation of weak solutions are weak solutions as well.

\begin{lemma}
\label{Concatenation}Let $f\in L^{2}(\tau,T;V^{\prime})$, $\beta\geq
1,\ \tau<s<T$.

\begin{itemize}
\item If $u$ is a weak solution on $[\tau,T]$, then $v\left(
\text{\textperiodcentered}\right)  =u\left(  \text{\textperiodcentered
}+s\right)  $ is a weak solution on $[\tau+s,T+s].$

\item If $u$ is a weak solution on $[\tau,s]$ and $v$ is a weak solution on
$[s,T]$ with $v\left(  s\right)  =u\left(  s\right)  $, then%
\[
w(t)=\left\{
\begin{array}
[c]{c}%
u\left(  t\right)  \text{ if }\tau\leq t\leq s,\\
v\left(  t\right)  \text{ if }s\leq t\leq T,
\end{array}
\right.
\]
is a weak solution on $[\tau,T]$.
\end{itemize}
\end{lemma}

If a function $u:[\tau,+\infty)\rightarrow H$ is a weak solution for any
$T>\tau$, then it is called a global solution. By Lemma \ref{Concatenation} if
$f\in L^{2}(\tau,T;V^{\prime})$, for all $T>\tau$, then any weak solution can
be extended to a global one \cite[p. 3580]{bib32}, so let%
\[
D_{\tau}(u_{\tau})=\{{u(\cdot):u(\cdot)}\text{ is a global weak solution to
(\ref{eq:NSDE-1-intro}) satisfying \eqref{eq:NSDE-4-Eineq} and }%
u(\tau)=u_{\tau}\}
\]
and denote the corresponding attainability set for $t\geq\tau$ by
\[
K_{t}(u_{\tau})={u(t):u(\cdot)\in D_{\tau}(u_{\tau})}.
\]

Let $Y=V^{\prime}+\left(  L^{\frac{\beta+1}{\beta}}(\Omega)\right)  ^{3}$
which is the dual space of $V\cap(L^{\beta+1}(\Omega))^{3}$. By standard
estimates on $B$ for any weak solution we have that%
\[
Au\in L^{2}(\tau,T;V^{\prime}),
\]%
\[
B(u,u)\in L^{4/3}(\tau,T;V^{\prime}),
\]%
\[
P\alpha|u|^{\beta-1}u=G(u)\in L^{\frac{\beta+1}{\beta}}(\tau,T;(L^{\frac
{\beta+1}{\beta}}(\Omega))^{3}),
\]
which implies in particular that
\[
-\mu Au-B(u,u)-G(u)+f\in L^{\frac{\beta+1}{\beta}}(\tau,T;(L^{\frac{\beta
+1}{\beta}}(\Omega))^{3})+L^{\frac{4}{3}}(\tau,T;V^{\prime})\subset L^{1}%
(\tau,T;Y).
\]
It follows from equality \eqref{eq:NSDE-3-weaksol} and a standard result
\cite[p.250 Lemma 1.1]{bib18} that%
\begin{equation}
\frac{du}{dt}=-\mu Au-B(u,u)-G(u)+f \label{eq:NSDE-5-weaksol}%
\end{equation}
in the sense of $Y$-valued distributions. Hence, the derivative $\frac{du}%
{dt}$ belongs to the space $Y$ and equality \eqref{eq:NSDE-5-weaksol} is
satisfied in the space $Y$ for a.a $t\in(0,T)$. In order to obtain good
estimates of weak solutions we need $\frac{du}{dt}$ to be more regular.

\begin{lemma}
\label{Deriv} Let $u$ be a weak solution to (\ref{eq:NSDE-1-intro}) such a
that $u\in L^{q}(\tau,T;(L^{q}(\Omega))^{3})$ which $q\geq4$. Then
\begin{equation}
\frac{du}{dt}\in L^{2}(\tau,T;V)+L^{\frac{\beta+1}{\beta}}(\tau,T;(L^{\frac
{\beta+1}{\beta}}(\Omega))^{3}),\ u\in C([0,T],H), \label{eq:NSDE-6}%
\end{equation}
the map $t\rightarrow|u(t)|^{2}$ is absolutely continuous and
\begin{equation}
\frac{d}{dt}|u(t)|^{2}=2\langle u,\frac{du}{dt}\rangle\text{ for a.a. }%
t\in(\tau,T). \label{eq:NSDE-7}%
\end{equation}

\end{lemma}

\begin{corollary}
\label{DerivCor}If $\beta\geq3$, any weak solution to (\ref{eq:NSDE-1-intro})
satisfies \eqref{eq:NSDE-6}-\eqref{eq:NSDE-7}.
\end{corollary}

\begin{corollary}
\label{EnergyEq}If $\beta\geq3$, any weak solution to (\ref{eq:NSDE-1-intro})
satisfies
\[
V_{\tau}(u(t))=V_{\tau}(u(s))\text{ for all }t\geq s\geq\tau.
\]

\end{corollary}

\begin{proof}
This is easily obtained by multiplying the equation by $u$, using
(\ref{eq:NSDE-7}) and integrating.
\end{proof}

\bigskip

If we assume $\beta>3,\alpha>0$ or $\beta=3,4\alpha\mu\geq1,$ for $u_{0}\in
H,f\in L^{2}(\tau,T;V^{\prime})$, then the weak solution to
(\ref{eq:NSDE-1-intro}) is unique \cite{bib26}. This unique weak solution is
continuous with respect to the initial datum $u_{0}$ (see \cite{bib27},
\cite{bib32}).

\section{The Globally Modified Navier-Stokes Equations with Damping (GMNSED)}

In this section, we introduce the Globally Modified Navier-Stokes Equations
with Damping (GMNSED) as an auxiliary tool in order to study the Kneser
property for problem (\ref{eq:NSDE-1-intro}) and study the existence,
uniqueness and continuity of solutions.

The Cauchy problem for the GMNSED is given by%
\begin{equation}%
\begin{cases}
u_{t}-\mu\Delta u+F_{N}(\left\Vert u\right\Vert )(u\cdot\nabla)u+\nabla
p+\alpha F_{N}(\left\Vert u\right\Vert ^{\beta-1})|u|^{\beta-1}u=f,\\
\operatorname{div}\text{ }u=0,\\
u|_{\partial\Omega}=0,\\
u(\tau)=u_{\tau},
\end{cases}
\label{eq:NSDE-8-GMNSED}%
\end{equation}
where $F_{N}:\mathbb{R}^{+}\rightarrow\mathbb{R}^{+}$ is defined by
$F_{N}(r):=min(1,\frac{N}{r}),r\in\mathbb{R}^{+}$. Here, $N$ is a fixed
positive number. They were introduced for the Navier-Stokes system in
\cite{bib28}. The difference from the Navier-Stokes system comes from the fact
that we multiply the term $\alpha|u|^{\beta-1}u$ by $F_{N}(\left\Vert
u\right\Vert ^{\beta-1})$. In this way, we kame both nonlinear terms more
regular and we are able to derive beter estimates for the solutions. Thus,
this model is a natural extension of the Globally Modified Navier-Stokes Equations.

A weak solution of the GMNSED \eqref{eq:NSDE-8-GMNSED} is defined as for
system (\ref{eq:NSDE-1-intro}), but with the operator $B(u,v)$ replaced by
$B_{N}(u,v)$, the element of $V^{\prime}$ defined by $\langle B_{N}%
(u,v),w\rangle=b_{N}(u,v,w)$, for all $u,v,w\in V$, where $b_{N}%
(u,v,w)=F_{N}(||v||)b(u,v,w)$. Also, we replace the operator $G(u)$ by
$G_{N}(u)=F_{N}(||u||^{\beta-1})G(u)$ and $u\in L^{\beta+1}(\tau
,T;(L^{\beta+1}(\Omega)^{3}))$ by $F_{N}(||u||^{\beta-1})u\in L^{\beta+1}%
(\tau,T;(L^{\beta+1}(\Omega)^{3}))$. It was shown in \cite{bib28} that
\begin{equation}
||B_{N}(u,v)||_{V^{\prime}}\leq NC_{1}||u||. \label{eq:NSDE-9}%
\end{equation}
Also, it is known \cite{bib37} that
\begin{equation}
|b(u,v,w)|\leq C_{3}||u||^{1/2}|Au|^{1/2}||v|||w|,\forall u\in D(A),v\in
V,w\in H. \label{eq:NSDE-9C}%
\end{equation}
It follows from (\ref{eq:NSDE-9}) that if $f\in L^{2}(\tau,T;V^{\prime})$,
then
\[
\frac{du}{dt}\in L^{2}(\tau,T;V^{\prime})+L^{\frac{\beta+1}{\beta}}%
(\tau,T;(L^{\frac{\beta+1}{\beta}}(\Omega)^{3})),
\]
so $u\in C(\tau,T;H)$ and (\ref{eq:NSDE-7}) hold \cite[Chapter II, Theorem
1.8]{bib03}. Hence, multiplying the equation by $u$ it follows that any weak
solution to (\ref{eq:NSDE-8-GMNSED}) satisfies the energy equality%
\begin{equation}
V_{\tau}^{N}(u(t))=V_{\tau}^{N}(u(s)),\ \forall\ \tau\leq s\leq t,
\label{eq:NSDE-N-Eineq}%
\end{equation}
where
\[
V_{\tau}^{N}(u(t)):=\frac{1}{2}\left\vert u(t)\right\vert ^{2}+\mu\int_{\tau
}^{t}||u(r)||^{2}+\alpha\int_{\tau}^{t}F_{N}(||u(r)||^{\beta-1})|u(r)|^{\beta
+1}dr-\int_{\tau}^{t}\langle f(r),u(r)\rangle dr.
\]

\begin{theorem}
\label{thm:Th03-primero}Let $1\leq\beta\leq3$. If $f\in L^{2}(\tau,T;H)$ for
all $T>\tau$ and if $u_{\tau}\in V$, then there exists a global strong
solution of system \eqref{eq:NSDE-8-GMNSED}, that is, a global weak solution
such that
\[
u\in L^{\infty}(\tau,T;V)\cap L^{2}(\tau,T;D(A))\text{ for all }T>\tau.
\]
Moreover, $u\in C([\tau,T],V)$ for all $T>\tau.$
\end{theorem}

\begin{proof}
Let $T>\tau$.

Following \cite{bib28} we use the Garlekin approximation given by
$u_{m}(x,t)=\sum_{j=1}^{m}u_{mj}(t)\phi_{j}(x)$ and $Au_{m}(x,t)=\sum
_{j=1}^{m}\lambda_{j}u_{mj}(t)\phi_{j}(x)$. Here the $\lambda_{j}$ and
$\phi_{j}$ are the corresponding eigenvalues and orthonormal eigenfunctions of
the operator $A$ and $P_{m}$ is the projection onto the subspace of $H$
spanned by {\ $\phi_{1},\phi_{2},...,\phi_{m}$ }. On the other hand, we have
$|u_{m}(x,t)|^{2}=\sum_{j=1}^{m}u_{mj}^{2}$, $|Au_{m}(x,t)|^{2}=\sum_{j=1}%
^{m}\lambda_{j}^{2}u_{mj}^{2}$ and $||u_{m}(x,t)||^{2}=\sum_{j=1}^{m}%
\lambda_{j}u_{mj}^{2}$. The expressions $|u_{m}|,||u_{m}||$ can be interpreted
as either the euclidean norm of $u_{m}\in\mathbb{R}^{m}$ or the $L^{2}$-norm
of $u_{m}\in H$ or $V$.

The functions $u_{m}$ satisfy the equality
\begin{equation}
\frac{d}{dt}u_{m}+\mu Au_{m}+P_{m}B_{N}(u_{m},u_{m})+P_{m}G_{N}(u_{m})=P_{m}f
\label{eq:NSDE-10}%
\end{equation}
with $G_{N}(u_{m})=F_{N}(||u_{m}||^{\beta-1})G(u_{m})$ and $u_{m}(\tau
)=P_{m}u_{\tau}$. If we take the inner product with $u_{m}$ and use that
$b(u_{m},u_{m},u_{m})=0$ we obtain that
\begin{equation}
\frac{d}{dt}|u_{m}|^{2}+\mu||u_{m}||^{2}+2\alpha F_{N}(||u_{m}||^{\beta
-1})\left\Vert u_{m}\right\Vert _{\beta+1}^{\beta+1}\leq\frac{|f|^{2}}%
{\mu\lambda_{1}} \label{eq:NSDE-11}%
\end{equation}
and
\[
\frac{d}{dt}|u_{m}|^{2}+\mu\lambda_{1}|u_{m}|^{2}\leq\frac{|f|^{2}}{\mu
\lambda_{1}}.
\]
Integrating \eqref{eq:NSDE-11} over $(\tau,t)$ we obtain
\begin{equation}
|u_{m}(t)|^{2}+\mu\int_{\tau}^{t}||u_{m}(s)||^{2}ds+2\alpha\int_{\tau}%
^{t}F_{N}(||u_{m}||^{\beta-1})|u_{m}(s)|_{\beta+1}^{\beta+1}ds\leq|u_{m\tau
}|^{2}+\int_{\tau}^{t}\frac{|f(s)|^{2}}{\lambda_{1}\mu}ds \label{eq:NSDE-12b}%
\end{equation}
Hence, there exists $u\in L^{2}(\tau,T;V)\cap L^{\infty}(\tau,T;H)$, and a
subsequence of $\{u_{m}\}$ which converges weak-star to $u$ in ${L^{\infty
}(0,T;H)}$, and weakly in $L^{2}(0,T;V)$. By the Aubin-Lions Compactness
Theorem \cite[Chapter 1]{bib11}, it follows that a subsequence in fact
converges strongly to $u$ in $L^{2}(0;T;H)$ and a.e. in $(\tau,T)\times\Omega
$. However, the weak convergence in $L^{2}(\tau,T,V)$ is not enough to ensure
that $||u_{m}(t)||\rightarrow||u(t)||$ or at least $F_{N}(||u_{m}%
(t)||)\rightarrow F_{N}(||u(t)||)$ and $F_{N}(||u_{m}(t)||^{\beta
-1})\rightarrow F_{N}(||u(t)||^{\beta-1})$ for a.a. $t$.

Thus we go one step further and find a stronger estimate. We now take the
inner product of the Galerkin ODE \eqref{eq:NSDE-10} with $Au_{m}$ and obtain
\begin{equation}
\frac{1}{2}\frac{d}{dt}||u_{m}||^{2}+\mu|Au_{m}|^{2}+b_{N}(u_{m},u_{m}%
,Au_{m})+\alpha F_{N}(||u_{m}||^{\beta-1})\int_{\Omega}|u_{m}|^{\beta-1}%
u_{m}Au_{m}dx=(f,Au_{m}). \label{eq:NSDE-13}%
\end{equation}
As $F_{N}(r)r\leq N$ we can estimate the term with damping as follows:%
\begin{equation}
\alpha F_{N}(||u_{m}||^{\beta-1})|\int_{\Omega}|u_{m}|^{\beta-1}u_{m}%
Au_{m}dx|\leq\frac{\alpha N}{||u_{m}||^{\beta-1}}\int_{\Omega}|u_{m}|^{\beta
}|Au_{m}|dx\leq\frac{{\alpha}^{2}N^{2}}{2\mu}\frac{|u_{m}|_{2\beta}^{2\beta}%
}{||u_{m}||^{2\beta-2}}+\frac{\mu}{2}|Au_{m}|^{2}. \label{eq:NSDE-13b}%
\end{equation}
As%
\[
\left\vert (f,Au_{m})\right\vert \leq\frac{\left\vert f\right\vert ^{2}}{\mu
}+\frac{\mu}{4}\left\vert Au_{m}\right\vert ^{2},
\]
using $F_{N}(r)r\leq N$, \eqref{eq:NSDE-9C} and the embedding $V\subset
(L^{2\beta}(\Omega))^{3}$ we have:
\[
\frac{1}{2}\frac{d}{dt}||u_{m}||^{2}+\mu|Au_{m}|^{2}\leq\frac{|f|^{2}}{\mu
}+\frac{7\mu}{8}|Au_{m}|^{2}+D_{N}||u_{m}||^{2}+C\frac{\alpha^{2}N^{2}}{2\mu
}||u_{m}||^{2},
\]%
\begin{equation}
\frac{1}{2}\frac{d}{dt}||u_{m}||^{2}+\frac{\mu}{8}|Au_{m}|^{2}\leq
\frac{|f|^{2}}{\mu}+C_{N}||u_{m}||^{2}. \label{eq:NSDE-14}%
\end{equation}
Finally:
\begin{equation}
||u_{m}(t)||^{2}+\frac{\mu}{4}\int_{\tau}^{t}|Au_{m}|^{2}ds\leq\frac{2}{\mu
}\int_{\tau}^{t}|f(s)|^{2}ds+2C_{N}\int_{\tau}^{t}||u_{m}||^{2}ds+||u_{m}%
(\tau)||^{2}. \label{eq:NSDE-14B}%
\end{equation}
Then, using Gronwall's lemma and the compactness theorem, up to a subsequence,
we get%
\[
u_{m}\rightarrow u\text{ strongly in }L^{2}(\tau,T;V),
\]%
\[
u_{m}\rightarrow u\text{ a.e. in }(\tau,T)\times\Omega,
\]%
\[
u_{m}\rightharpoonup u\text{ weakly in }L^{2}(\tau,T;D(A)),
\]%
\[
u_{m}\rightharpoonup u\text{ weakly-star in }L^{\infty}(\tau,T;V),
\]%
\[
\frac{du_{m}}{dt}\rightharpoonup\frac{du}{dt}\text{ weakly in }L^{2}%
(\tau,T;H),
\]
This last expression is obtained by observing that $|G_{N}(u_{m})|\leq
\alpha|u_{m}|_{2\beta}^{\beta}\leq\alpha C||u||^{\beta}$, and $|B_{N}%
(u_{m},u_{m})|\leq C||u_{m}||^{3/2}|Au_{m}|^{1/2}$, where we have used \eqref{eq:NSDE-9C}.

The proof of the fact that the limit function $u$ is a weak solution
satisfying the initial condition is similar to the classical one for the
Navier-Stokes system (see \cite{bib11}, \cite{bib18}). We only need to pay
attention to the passage to the limit for the terms $P_{m}B_{N}(u_{m},u_{m})$
and $P_{m}G_{N}(u_{m}).$

First, for the term $P_{m}B_{N}(u_{m},u_{m})$, by linearity, density and the
properties of $b,$ it is enough to prove that%
\begin{equation}
\lim_{m\rightarrow\infty}\int_{\tau}^{t}F_{N}(||u_{m}(s)||)b(u_{m}(s),\phi
_{j},u_{m}(s))ds=\int_{\tau}^{t}F_{N}(||u(s)||)b(u(s),\phi_{j},u(s))ds
\label{LimBN}%
\end{equation}
for $j\geq1,t\geq\tau$. Let $T>\tau$ be fixed. Since $u_{m}$ is bounded in
$L^{\infty}(\tau,T;H)\cap L^{2}(\tau,T;V),$ according to \cite{bib11} for each
$1\leq i$, $k\leq3$, the product of the components $u_{m}^{(i)},u_{m}^{(k)}$
is bounded in $L^{2}(\tau,T;L^{3/2}(\Omega))$, and, therefore, as $0\leq
F_{N}\left(  r\right)  \leq1$, the product $F_{N}(||u_{m}($\textperiodcentered
$)||)u_{m}^{(i)}u_{m}^{(k)}$ is also bounded in $L^{2}(\tau,T;(L^{3/2}%
(\Omega)))$. But then, up to a subsequence, we can assume that%
\[
F_{N}(||u_{m}(\text{\textperiodcentered})||)u_{m}^{(i)}u_{m}^{(k)}%
\rightarrow\chi_{i}{}_{k}\text{ weak in }L^{2}(\tau,T;(L^{3/2}(\Omega
)))\text{, }1\leq i,k\leq3.
\]
Since $u_{m}\rightarrow u$ a.e. in $(\tau,T)\times\Omega$ and $F_{N}%
(||u_{m}(s)||)\rightarrow F_{N}(||u(s)||)$ for a.a. $s\in\left(
\tau,T\right)  $, we infer that $F_{N}(||u_{m}(s)||)u_{m}^{(i)}u_{m}%
^{(k)}\rightarrow F_{N}(||u(s)||)u^{(i)}u^{(k)}$, a.e. in $(\tau
,T)\times\Omega$, for $1\leq i$, $k\leq3$. The sequence $F_{N}(||u_{m}%
($\textperiodcentered$)||)u_{m}^{(i)}u_{m}^{(k)}$ is also bounded in
$L^{3/2}(\tau,T)\times\Omega)$, and then $F_{N}(||u_{m}($\textperiodcentered
$)||)u_{m}^{(i)}u_{m}^{(k)}\rightarrow F_{N}(||u(s)||)u^{(i)}u^{(k)}$ weakly
in $L^{3/2}(\tau,T)\times\Omega),$ $1\leq i$, $k\leq3$ \cite[Lemma 1.3]%
{bib11}. This implies that $\chi_{i}{}_{k}=F_{N}(||u($\textperiodcentered
$)||)u^{(i)}u^{(k)}$. Hence, (\ref{LimBN}) is obtained proceeding as in
\cite{bib11} for the Navier-Stokes system.

Second, for the term $P_{m}G_{N}(u_{m})$ it is enough to prove that
\begin{equation}
P_{m}G_{N}(u_{m})\rightarrow G_{N}(u)\text{ weakly in }L^{q}(\tau,T;\left(
L^{q}(\Omega)\right)  ^{3}), \label{ConvergGN}%
\end{equation}
where $q=(\beta+1)^{\prime}=\frac{\beta+1}{\beta}$. We know from
\eqref{eq:NSDE-12b} that
\[
\int_{\tau}^{T}F_{N}(||u_{m}(s)||^{\beta-1})|u_{m}(s)|_{\beta+1}^{\beta
+1}ds\leq C,
\]
so using $0\leq F_{N}\left(  r\right)  \leq1$ we get
\begin{align*}
\int_{\tau}^{T}|F_{N}(||u_{m}(s)||^{\beta-1})u_{m}(s)|_{\beta+1}^{\beta+1}ds
&  \leq\int_{\tau}^{T}F_{N}(||u_{m}(s)||^{\beta-1})|u_{m}(s)|_{\beta+1}%
^{\beta+1}ds\leq C,\\
\int_{\tau}^{T}|F_{N}(||u_{m}(s)||^{\beta-1})|u_{m}(s)|^{\beta-1}u_{m}%
(s)|_{q}^{q}ds  &  \leq\int_{\tau}^{T}F_{N}(||u_{m}(s)||^{\beta-1}%
)|u_{m}(s)|_{\beta+1}^{\beta+1}ds\leq C,
\end{align*}
and then%
\begin{equation}
F_{N}(||u_{m}||^{\beta-1})u_{m}\text{ is bounded in }L^{\beta+1}%
(\tau,T;\left(  L^{\beta+1}(\Omega)\right)  ^{3}), \label{FNumBounded}%
\end{equation}%
\begin{equation}
F_{N}(||u_{m}||^{\beta-1})|u_{m}|^{\beta-1}u_{m}\text{ is bounded in }%
L^{q}(\tau,T;\left(  L^{q}(\Omega)\right)  ^{3}), \label{GmBoundedLq}%
\end{equation}%
\[
F_{N}(||u_{m}||^{\beta-1})u_{m}\rightarrow\chi_{1}\text{ weakly in }%
L^{\beta+1}(\tau,T;\left(  L^{\beta+1}(\Omega)\right)  ^{3}),
\]%
\[
F_{N}(||u_{m}||^{\beta-1})|u_{m}|^{\beta-1}u_{m}\rightarrow\chi_{2}\text{
weakly in }L^{q}(\tau,T;\left(  L^{q}(\Omega)\right)  ^{3}),
\]
for some $\chi_{1},\chi_{2}$. We need to show that%
\begin{equation}
\chi_{1}=F_{N}(||u||^{\beta-1})u,\ \chi_{2}=F_{N}(||u||^{\beta-1}%
)|u|^{\beta-1}u. \label{eq:24-v}%
\end{equation}
Since $u_{m}\rightarrow u$ a.e. in $(\tau,T)\times\Omega$ and $F_{N}%
(||u_{m}(s)||)\rightarrow F_{N}(||u(s)||)$ for a.a. $s\in\left(
\tau,T\right)  $, we deduce that
\begin{align*}
F_{N}(||u_{m}||^{\beta-1})||u_{m}|^{\beta-1}u_{m}  &  \rightarrow
F_{N}(||u||^{\beta-1})|u|^{\beta-1}u,\\
\ F_{N}(||u_{m}||^{\beta-1})u_{m}  &  \rightarrow F_{N}(||u||^{\beta
-1})u\text{ a.e. in }(\tau,T)\times\Omega.
\end{align*}
By (\ref{FNumBounded}), (\ref{GmBoundedLq}), equality \eqref{eq:24-v} follows
from Lemma 1.3 in \cite{bib11}. We have obtained that%
\[
G_{N}(u_{m})\rightarrow G_{N}(u)\text{ weakly in }L^{q}(\tau,T;\left(
L^{q}(\Omega)\right)  ^{3}),
\]%
\[
F_{N}(\left\Vert u_{m}\right\Vert ^{\beta-1})u_{m}\rightarrow F_{N}(\left\Vert
u\right\Vert ^{\beta-1})u\text{ weakly in }L^{\beta+1}(\tau,T;\left(
L^{\beta+1}(\Omega)\right)  ^{3}).
\]
The convergence (\ref{ConvergGN}) is established in a similar way as in
\cite[Theorem 8.4]{bib14}.

Hence, we have proved that $u$ is a strong solution on $[\tau,T]$. Finally,
since $u\in L^{2}(\tau,T;D(A))$ and $\dfrac{du}{dt}\in L^{2}(\tau,T;H)$, a
standard result \cite{bib18} implies that $u\in C([\tau,T],V).$

Finally, by a diagonal argument we extend the solution to $[0,+\infty).$
\end{proof}

\bigskip

\begin{theorem}
\label{thm:Th04-segundo} Let $1\leq\beta\leq3$. If $f\in L^{2}(0,T;H)$ for all
$\tau<T$ and, if $u_{\tau}\in H$, there exists a global weak solution which is
a strong one in the following sense:%
\[
u\in L^{\infty}(\tau+\epsilon,T;V)\cap L^{2}(\tau+\epsilon,T;D(A)),
\]
for any $\epsilon\in(0,T-\tau),\ T>\tau.$
\end{theorem}

\begin{proof}
Fix $T>\tau$. For $\tau\leq t\leq T-r,\ r>0$, integrating over $\left(
t,t+r\right)  $ in (\ref{eq:NSDE-11}) we have%
\[
\mu\int_{t}^{t+r}\left\Vert u_{m}(s)\right\Vert ^{2}ds\leq\left\vert
u_{m}(t)\right\vert ^{2}+\frac{1}{\mu\lambda_{1}}\int_{t}^{t+r}|f(s)|^{2}%
ds\leq\left\vert u_{m}(t)\right\vert ^{2}+\frac{1}{\mu\lambda_{1}}\int_{\tau
}^{T}|f(s)|^{2}ds.
\]
By (\ref{eq:NSDE-12b}) we obtain%
\[
\mu\int_{t}^{t+r}\left\Vert u_{m}(s)\right\Vert ^{2}ds\leq|u_{m\tau}%
|^{2}+\frac{1}{\mu\lambda_{1}}\int_{\tau}^{T}|f(s)|^{2}ds=|u_{m\tau}%
|^{2}+C_{T}.
\]
It follows the existence of $t_{0}\in(t,t+r)$ such that
\[
||u_{m}(t_{0})||^{2}\leq\frac{2}{\mu r}(|u_{m\tau}|^{2}+C_{T}).
\]
Thus we can ensure also that for any $\tau+\epsilon\leq t\leq T$, $\epsilon>0$
there exists a $t_{0}\in(t-\epsilon,t)$ such that%
\begin{equation}
||u_{m}(t_{0})||^{2}\leq\frac{2}{\mu\epsilon}(|u_{m\tau}|^{2}+C_{T}).
\label{eq:NEW19v}%
\end{equation}
From (\ref{eq:NSDE-14}) we have%
\begin{align*}
||u_{m}(t)||^{2}  &  \leq e^{2C_{N}(t-t_{0})}||u_{m}(t_{0})||^{2}+\frac{2}%
{\mu}\int_{t_{0}}^{t}e^{2C_{N}(t-s)}|f(s)|^{2}ds\\
&  \leq e^{2\epsilon C_{N}}(\frac{2}{\mu\epsilon}(|u_{m\tau}|^{2}+C_{T}%
)+\frac{2}{\mu}\int_{\tau}^{T}|f(s)|^{2}ds)\text{ for }\tau+\epsilon\leq t\leq
T.
\end{align*}
The last inequality, together with (\ref{eq:NSDE-12b}) and (\ref{eq:NSDE-14})
implies that the sequence ${u_{m}}$ is bounded in $L^{\infty}(\tau,T;H)$,
$L^{2}(\tau,T;V)$, $L^{\infty}(\tau+\epsilon,T;V)$ and $L^{2}(\tau
+\epsilon;T;D(A))$. Then, arguing as in Theorem \ref{thm:Th03-primero} we see
that the sequence $\left\{  \dfrac{du_{m}}{dt}\right\}  $ is also bounded in
$L^{2}(\tau+\epsilon;T;H)$. Hence, there exist an element $u\in L^{\infty
}(\tau,T,H)\cap L^{2}(\tau,T,V)\cap L^{\infty}(\tau+\epsilon,T,V)\cap
L^{2}(\tau+\epsilon,T,D(A))$ and a subsequence of $\{{u_{m}\}}$, that we will
also denote by $\{u_{m}\}$, such that%
\[
u_{m}\rightharpoonup u\text{ weakly in }L^{2}(\tau,T;V),
\]%
\[
u_{m}\rightharpoonup u\text{ weakly-star in }L^{\infty}(\tau,T;H),
\]%
\[
u_{m}\rightarrow u\text{ strongly in }L^{2}(\tau,T;H),
\]%
\[
u_{m}(t)\rightarrow u(t)\text{ strongly in }H\text{ for a.a. }t\in\left(
\tau,T\right)  ,
\]%
\[
u_{m}\rightarrow u\text{ a.e. in }(\tau,T)\times\Omega,
\]%
\[
u_{m}\rightarrow u\text{ strongly in }L^{2}(\tau+\epsilon,T;V),
\]%
\[
u_{m}\rightharpoonup u\text{ weakly in }L^{2}(\tau+\epsilon,T;D(A)),
\]%
\[
u_{m}\rightharpoonup u\text{ weakly-star in }L^{\infty}(\tau+\epsilon,T;V),
\]%
\[
\frac{du_{m}}{dt}\rightharpoonup\frac{du}{dt}\text{ weakly in }L^{2}%
(\tau+\epsilon,T;H)
\]

The proof of the fact that $u$ is a weak solution on $[\tau,T]$ and its
regularity follows the same lines as in Theorem \ref{thm:Th03-primero}. By a
diagonal argument we extend the result for any $T>\tau.$
\end{proof}

\bigskip

Now, we have to show that the weak solution is unique.

\begin{lemma}
\label{FNLemma}Let $\beta\geq1$. For any $u,v\in V$ we have:

\begin{enumerate}
\item If $\left\Vert u\right\Vert ^{\beta-1},\left\Vert v\right\Vert
^{\beta-1}\leq N$, then
\[
F_{N}(\left\Vert u\right\Vert ^{\beta-1})-F_{N}(\left\Vert v\right\Vert
^{\beta-1})=0.
\]

\item If $\left\Vert u\right\Vert ^{\beta-1}\geq\left\Vert v\right\Vert
^{\beta-1}>N$, then
\[
\left\vert F_{N}(\left\Vert u\right\Vert ^{\beta-1})-F_{N}(\left\Vert
v\right\Vert ^{\beta-1})\right\vert \leq C(\beta)\frac{2N}{\left\Vert
v\right\Vert ^{\beta}} \left\Vert u-v\right\Vert \leq C(\beta)\frac
{1}{\left\Vert v\right\Vert } \left\Vert u-v\right\Vert .
\]

\item If $\left\Vert u\right\Vert ^{\beta-1}>N,\left\Vert v\right\Vert
^{\beta-1}\leq N$, then%
\[
\left\vert F_{N}(\left\Vert u\right\Vert ^{\beta-1})-F_{N}(\left\Vert
v\right\Vert ^{\beta-1})\right\vert \leq C(\beta)\frac{1}{\left\Vert
u\right\Vert }\left\Vert u-v\right\Vert ,
\]

\end{enumerate}
\end{lemma}

\begin{proof}
If $\left\Vert u\right\Vert ^{\beta-1},\left\Vert v\right\Vert ^{\beta-1}\leq
N$, then $F_{N}(\left\Vert u\right\Vert ^{\beta-1})=F_{N}(\left\Vert
v\right\Vert ^{\beta-1})=1$, so $F_{N}(\left\Vert u\right\Vert ^{\beta
-1})-F_{N}(\left\Vert v\right\Vert ^{\beta-1})=0.$

If $\left\Vert u\right\Vert ^{\beta-1}\geq\left\Vert v\right\Vert ^{\beta
-1}>N$, then%
\begin{align*}
F_{N}(\left\Vert u\right\Vert ^{\beta-1})-F_{N}(\left\Vert v\right\Vert
^{\beta-1})  &  =N\left(  \frac{1}{\left\Vert u\right\Vert ^{\beta-1}}%
-\frac{1}{\left\Vert v\right\Vert ^{\beta-1}}\right) \\
&  =\frac{N}{\left\Vert u\right\Vert ^{\beta-1}\left\Vert v\right\Vert
^{\beta-1}}\left(  \left\Vert v\right\Vert ^{\beta-1}-\left\Vert u\right\Vert
^{\beta-1}\right)  .
\end{align*}
There is $\theta\in(0,1)$ such that%
\[
\left\Vert v\right\Vert ^{\beta-1}-\left\Vert u\right\Vert ^{\beta-1}=\left(
\beta-1\right)  \left(  \theta\left\Vert u\right\Vert +(1-\theta)\left\Vert
v\right\Vert \right)  ^{\beta-2}\left(  \left\Vert u\right\Vert -\left\Vert
v\right\Vert \right)  .
\]
Then
\begin{align*}
\left\vert F_{N}(\left\Vert u\right\Vert ^{\beta-1})-F_{N}(\left\Vert
v\right\Vert ^{\beta-1})\right\vert  &  \leq C(\beta)\frac{N}{\left\Vert
u\right\Vert ^{\beta-1}\left\Vert v\right\Vert ^{\beta-1}}\left(  \left\Vert
u\right\Vert ^{\beta-2}+\left\Vert v\right\Vert ^{\beta-2}\right)  \left\Vert
u-v\right\Vert \\
&  =C(\beta)N\left(  \frac{1}{\left\Vert u\right\Vert \left\Vert v\right\Vert
^{\beta-1}}+\frac{1}{\left\Vert v\right\Vert \left\Vert u\right\Vert
^{\beta-1}}\right)  \left\Vert u-v\right\Vert \\
&  \leq C(\beta)\frac{2N}{\left\Vert v\right\Vert ^{\beta}}\left\Vert
u-v\right\Vert \\
&  \leq C(\beta)\frac{2||v||^{\beta-1}}{\left\Vert v\right\Vert ^{\beta}%
}\left\Vert u-v\right\Vert \\
&  \leq C(\beta)\frac{2}{\left\Vert v\right\Vert }\left\Vert u-v\right\Vert .
\end{align*}

Finally, let $\left\Vert u\right\Vert ^{\beta-1}>N,\ \left\Vert v\right\Vert
^{\beta-1}\leq N.$ Hence,%

\begin{align*}
\left\vert F_{N}(\left\Vert u\right\Vert ^{\beta-1})-F_{N}(\left\Vert
v\right\Vert ^{\beta-1})\right\vert  &  =1-\frac{N}{\left\Vert u\right\Vert
^{\beta-1}}=\frac{1}{\left\Vert u\right\Vert ^{\beta-1}}\left(  \left\Vert
u\right\Vert ^{\beta-1}-N\right) \\
&  \leq\frac{1}{\left\Vert u\right\Vert ^{\beta-1}}\left(  \left\Vert
u\right\Vert ^{\beta-1}-\left\Vert v\right\Vert ^{\beta-1}\right) \\
&  =\frac{\left(  \beta-1\right)  }{\left\Vert u\right\Vert ^{\beta-1}}\left(
\theta\left\Vert u\right\Vert +(1-\theta)\left\Vert v\right\Vert \right)
^{\beta-2}\left(  \left\Vert u\right\Vert -\left\Vert v\right\Vert \right) \\
&  \leq C(\beta)\frac{\left\Vert u\right\Vert ^{\beta-2}+\left\Vert
v\right\Vert ^{\beta-2}}{\left\Vert u\right\Vert ^{\beta-1}}\left\Vert
u-v\right\Vert \\
&  \leq C(\beta)\frac{2}{\left\Vert u\right\Vert }\left\Vert u-v\right\Vert .
\end{align*}

\end{proof}

\begin{lemma}
\label{Lemma3} Let $1\leq\beta\leq3$. For any $\varepsilon>0$ there exists
$C\geq0$ (depending on $\beta,\alpha$ and $\varepsilon$) such that
\[
\alpha\int_{\Omega}\left(  F_{N}(\left\Vert u\right\Vert ^{\beta-1})\left\vert
u\right\vert ^{\beta-1}u-F_{N}(\left\Vert v\right\Vert ^{\beta-1})\left\vert
v\right\vert ^{\beta-1}v\right)  (u-v)dx\geq-\varepsilon\left\Vert
u-v\right\Vert ^{2}-CN^{2}\left\vert u-v\right\vert ^{2},
\]
for any $u,v\in V.$
\end{lemma}

\begin{proof}
For instance, let $\left\Vert u\right\Vert \geq\left\Vert v\right\Vert $.
First,
\begin{align*}
&  \int_{\Omega}\left(  F_{N}(\left\Vert u\right\Vert ^{\beta-1})\left\vert
u\right\vert ^{\beta-1}u-F_{N}(\left\Vert v\right\Vert ^{\beta-1})\left\vert
v\right\vert ^{\beta-1}v\right)  (u-v)dx\\
&  =\int_{\Omega}F_{N}(\left\Vert u\right\Vert ^{\beta-1})\left(  \left\vert
u\right\vert ^{\beta-1}u-\left\vert v\right\vert ^{\beta-1}v\right)  (u-v)dx\\
&  +\int_{\Omega}\left(  F_{N}(\left\Vert u\right\Vert ^{\beta-1}%
)-F_{N}(\left\Vert v\right\Vert ^{\beta-1}))\left\vert v\right\vert ^{\beta
-1}v\right)  (u-v)dx\\
&  \geq\int_{\Omega}\left(  F_{N}(\left\Vert u\right\Vert ^{\beta-1}%
)-F_{N}(\left\Vert v\right\Vert ^{\beta-1}))\left\vert v\right\vert ^{\beta
-1}v\right)  (u-v)dx.
\end{align*}
Second, we estimate the last integral in the previous inequality. If
$\left\Vert u\right\Vert ^{\beta-1},\left\Vert v\right\Vert ^{\beta-1}\leq N$,
then $F_{N}(\left\Vert u\right\Vert ^{\beta-1})-F_{N}(\left\Vert v\right\Vert
^{\beta-1})=0$, so the integral is equal to $0.$ Using Lemma \ref{FNLemma} and
the embedding $V\subset L^{2\beta}(\Omega)$, for $\left\Vert u\right\Vert
^{\beta-1},\left\Vert v\right\Vert ^{\beta-1}>N$ we have%
\begin{align*}
&  \left\vert \alpha\int_{\Omega}\left(  F_{N}(\left\Vert u\right\Vert
^{\beta-1})-F_{N}(\left\Vert v\right\Vert ^{\beta-1}))\left\vert v\right\vert
^{\beta-1}v\right)  (u-v)dx\right\vert \\
&  \leq\alpha\left\vert F_{N}(\left\Vert u\right\Vert ^{\beta-1}%
)-F_{N}(\left\Vert v\right\Vert ^{\beta-1})\right\vert \left\Vert v\right\Vert
_{2\beta}^{\beta}\left\vert u-v\right\vert \\
&  \leq\alpha C(\beta)\frac{2N}{\left\Vert v\right\Vert ^{\beta}}\left\Vert
u-v\right\Vert \left\Vert v\right\Vert _{2\beta}^{\beta}\left\vert
u-v\right\vert \\
&  \leq2\alpha\tilde{C}(\beta)N\left\Vert u-v\right\Vert \left\vert
u-v\right\vert \\
&  \leq\varepsilon\left\Vert u-v\right\Vert ^{2}+C(\beta,\varepsilon
,\alpha)N^{2}\left\vert u-v\right\vert ^{2}.
\end{align*}

Finally, if $\left\Vert u\right\Vert ^{\beta-1}>N,$ $||v||>0$ ( if $v=0$, the
result is obvious), $\left\Vert v\right\Vert ^{\beta-1}\leq N,$ then by Lemma
\ref{FNLemma} and the embedding $V\subset L^{2\beta}(\Omega)$ we get%
\begin{align*}
&  \left\vert \alpha\int_{\Omega}\left(  F_{N}(\left\Vert u\right\Vert
^{\beta-1})-F_{N}(\left\Vert v\right\Vert ^{\beta-1}))\left\vert v\right\vert
^{\beta-1}v\right)  (u-v)dx\right\vert \\
&  \leq\alpha\left\vert F_{N}(\left\Vert u\right\Vert ^{\beta-1}%
)-F_{N}(\left\Vert v\right\Vert ^{\beta-1})\right\vert \left\Vert v\right\Vert
_{2\beta}^{\beta}\left\vert u-v\right\vert \\
&  \leq\alpha C(\beta)\frac{1}{\left\Vert v\right\Vert }\left\Vert
u-v\right\Vert \left\Vert v\right\Vert ^{\beta}\left\vert u-v\right\vert \\
&  \leq2\alpha C(\beta)\left\Vert v\right\Vert ^{\beta-1}\left\Vert
u-v\right\Vert \left\vert u-v\right\vert \\
&  \leq\alpha C(\beta)N\left\Vert u-v\right\Vert \left\vert u-v\right\vert \\
&  \leq\varepsilon\left\Vert u-v\right\Vert ^{2}+C(\beta,\mu,\alpha
,\varepsilon)N^{2}\left\vert u-v\right\vert ^{2}.
\end{align*}

\end{proof}

\bigskip

\begin{theorem}
\label{thm:Th03-tercero} Let $1\leq\beta\leq3$. If $f\in L^{2}(\tau,T;H)$ for
all $\tau<T$ and, if $u_{\tau}\in H$, then the weak solution to
\eqref{eq:NSDE-8-GMNSED} is unique.
\end{theorem}

\begin{proof}
Let $u,v$ be two solutions; for $w=u-v$ from Lemma \ref{Lemma3} with
$\varepsilon=\frac{\mu}{2}$ we have%
\begin{align*}
&  \frac{1}{2}\frac{d}{dt}\left\vert w\right\vert ^{2}+\mu\left\Vert
w\right\Vert ^{2}+\left\langle F_{N}(\left\Vert u\right\Vert )\left(
u\text{\textperiodcentered}\nabla\right)  u-F_{N}(\left\Vert v\right\Vert
)\left(  v\text{\textperiodcentered}\nabla\right)  v,w\right\rangle \\
&  =-\alpha\int_{\Omega}\left(  F_{N}(\left\Vert u\right\Vert ^{\beta
-1})\left\vert u\right\vert ^{\beta-1}u-F_{N}(\left\Vert v\right\Vert
^{\beta-1})\left\vert v\right\vert ^{\beta-1}v\right)  (u-v)dx\\
&  \leq\frac{\mu}{2}\left\Vert w\right\Vert ^{2}+CN^{2}\left\vert w\right\vert
^{2}.
\end{align*}
We will estimate the term $\left\langle F_{N}(\left\Vert u\right\Vert )\left(
u\text{\textperiodcentered}\nabla\right)  u-F_{N}(\left\Vert v\right\Vert
)\left(  v\text{\textperiodcentered}\nabla\right)  v,w\right\rangle $ as in
\cite{bib29}. Using the inequality
\begin{equation}
|\langle B(u,v),w\rangle|\leq c_{0}\left\Vert u\right\Vert \left\Vert
v\right\Vert {\left\Vert w\right\Vert }^{1/2}{\left\vert w\right\vert }^{1/2},
\label{IneqBRomito}%
\end{equation}
see \cite{bib29}, and with an algebraic manipulation of $\mathcal{LN}%
=F_{N}(\left\Vert u\right\Vert )B(u,u)-F_{N}(\left\Vert v\right\Vert )B(v,v)$
we obtain
\[
\left\langle F_{N}(\left\Vert u\right\Vert )\left(  u\text{\textperiodcentered
}\nabla\right)  u-F_{N}(\left\Vert v\right\Vert )\left(
v\text{\textperiodcentered}\nabla\right)  v,u-v\right\rangle \leq\frac{\mu}%
{2}{\left\Vert {w}\right\Vert }^{2}+CN^{4}|w|^{2}.
\]

Applying Gronwall's lemma we obtain $w(t)=0$ for all $t\in\lbrack\tau,T]$, so
$u=v$, and the weak solution is unique.
\end{proof}

\bigskip

We finish this section by proving that the solutions to
(\ref{eq:NSDE-8-GMNSED}) depend continuously on $N$, $u_{\tau}$ and $f.$

\begin{lemma}
\label{IneqFMN}For any $M,N,p,r>0$ we have%
\[
\left\vert F_{M}(p)-F_{N}(r)\right\vert \leq\frac{M}{rp}\left\vert
p-r\right\vert +\frac{1}{r}\left\vert M-N\right\vert .
\]

\end{lemma}

\begin{proof}
We have
\[
\left\vert F_{M}(p)-F_{N}(r)\right\vert \leq\left\vert F_{M}(p)-F_{M}%
(r)\right\vert +\left\vert F_{M}(r)-F_{N}(r)\right\vert .
\]
For the first term, if $p,r\leq M$, then $F_{M}(p)-F_{M}(r)=0$. If $p,r>M$,
then
\[
\left\vert F_{M}(p)-F_{M}(r)\right\vert =M\left\vert \frac{1}{p}-\frac{1}%
{r}\right\vert =\frac{M}{rp}\left\vert p-r\right\vert .
\]
If, for instance, $p>M\geq r$, then%
\begin{align*}
\left\vert F_{M}(p)-F_{M}(r)\right\vert  &  =\left\vert \frac{M}%
{p}-1\right\vert =1-\frac{M}{p}=\frac{M}{rp}(\frac{pr}{M}-r)\\
&  \leq\frac{M}{rp}(p-r)=\frac{M}{rp}\left\vert p-r\right\vert .
\end{align*}
For the second term, the inequality
\[
\left\vert F_{M}(r)-F_{N}(r)\right\vert \leq\frac{1}{r}\left\vert
M-N\right\vert
\]
is proved in \cite[p.417]{bib28}.
\end{proof}

\begin{theorem}
\label{thm:Th03B} Let $1\leq\beta\leq3$. The weak solutions to
\eqref{eq:NSDE-8-GMNSED} depend continuously on $N$, $u_{\tau}$ and $f$. More
precisely, if $u,v$ are weak solutions corresponding, respectively, to
$N,u_{\tau},f_{u}$ and $M,v_{\tau},f_{v}$, then there exist constants
$C,\gamma(M,N)>0$ such that\
\[
|v\left(  t\right)  -u\left(  t\right)  |^{2}\leq Ce^{\gamma\left(
t-\tau\right)  }\left(  |v_{\tau}-u_{\tau}|^{2}+\left(  \int_{\tau}%
^{t}\left\Vert v\right\Vert ^{2}ds\right)  |M-N|^{2}+\int_{\tau}^{t}%
|f_{v}-f_{u}|^{2}ds\right)  \text{ for all }t\geq\tau.
\]

\end{theorem}

\begin{proof}
Considering two weak solutions $v,u$ to \eqref{eq:NSDE-8-GMNSED}, $w=v-u$
satisfies%
\[
\frac{d}{dt}|w|^{2}+2\mu||w||^{2}+2\left\langle F_{M}(\left\Vert v\right\Vert
)\left(  v\text{\textperiodcentered}\nabla\right)  v-F_{N}(\left\Vert
u\right\Vert )\left(  u\text{\textperiodcentered}\nabla\right)
u,w\right\rangle
\]%
\[
+2\alpha\int_{\Omega}\Big[F_{M}(||v||^{\beta-1})|v|^{\beta-1}v-F_{N}%
(||u||^{\beta-1})|u|^{\beta-1}u\Big]wdx=2(f_{v}-f_{u},w).
\]
Then%
\[
I=2\alpha F_{M}(||v||^{\beta-1})|v|^{\beta-1}v-F_{N}(||u||^{\beta
-1})|u|^{\beta-1}u
\]%
\[
=2\alpha F_{M}(||v||^{\beta-1})(|v|^{\beta-1}v-|u|^{\beta-1}u)+2\alpha
\big[F_{M}(||v||^{\beta-1})-F_{N}(||u||^{\beta-1})\big]|u|^{\beta-1}%
u=I_{1}+I_{2}.
\]
For $I_{1}$ we have%
\[
\int_{\Omega}I_{1}(v-u)dx=2\alpha\int_{\Omega}F_{M}(||v||^{\beta
-1})(|v|^{\beta-1}v-|u|^{\beta-1}u)(v-u)\geq0.
\]
For $I_{2}$ we obtain%
\[
\int_{\Omega}I_{2}(v-u)dx=2\alpha\int_{\Omega}\big[[F_{M}(||v||^{\beta
-1})-F_{M}(||u||^{\beta-1})]+[F_{M}(||u||^{\beta-1})-F_{N}(||u||^{\beta
-1})]\big]|u|^{\beta-1}u(v-u)dx
\]%
\[
=2\alpha\int_{\Omega}[F_{M}(||v||^{\beta-1})-F_{M}(||u||^{\beta-1}%
)]|u|^{\beta-1}u(v-u)dx+2\alpha\int_{\Omega}[F_{M}(||u||^{\beta-1}%
)-F_{N}(||u||^{\beta-1})]|u|^{\beta-1}u(v-u)dx.
\]
Using Lemma \ref{Lemma3} we have%
\[
2\alpha\int_{\Omega}[F_{M}(||v||^{\beta-1})-F_{M}(||u||^{\beta-1}%
)]|u|^{\beta-1}u(v-u)dxs\geq-\frac{\mu}{2}\left\Vert v-u\right\Vert ^{2}%
-C_{1}M^{2}\left\vert v-u\right\vert ^{2}.
\]
We consider the second integral. We take $u$ such that $\left\Vert
u\right\Vert >0$ (if $u=0$, then the integral is equal to $0$). Using Lemma 6
in \cite{bib28} we obtain%
\begin{equation}
\left\vert F_{M}(\left\Vert u\right\Vert ^{\beta-1})-F_{N}(\left\Vert
u\right\Vert ^{\beta-1})\right\vert \leq\frac{\left\vert M-N\right\vert
}{\left\Vert u\right\Vert ^{\beta-1}}. \label{IneqMN}%
\end{equation}
Thus,
\begin{align*}
&  \left\vert 2\alpha\int_{\Omega}\left(  F_{M}(\left\Vert u\right\Vert
^{\beta-1})-F_{N}(\left\Vert u\right\Vert ^{\beta-1})\right)  \left\vert
u\right\vert ^{\beta-1}u\left(  v-u\right)  dx\right\vert \\
&  \leq2\alpha\frac{\left\vert M-N\right\vert }{\left\Vert u\right\Vert
^{\beta-1}}\left\vert u\right\vert _{2\beta}^{\beta}\left\vert w\right\vert
\leq C_{2}\left\vert M-N\right\vert \left\Vert u\right\Vert \left\vert
w\right\vert \leq C_{3}\left\vert M-N\right\vert ^{2}\left\Vert u\right\Vert
^{2}+\left\vert w\right\vert ^{2}.
\end{align*}
Finally,%
\begin{align*}
&  2\alpha\int_{\Omega}\left(  F_{M}(\left\Vert v\right\Vert ^{\beta
-1})\left\vert v\right\vert ^{\beta-1}v-F_{N}(\left\Vert u\right\Vert
^{\beta-1})\left\vert u\right\vert ^{\beta-1}u\right)  \left(  v-u\right)
dx\\
&  \geq-\frac{\mu}{2}\left\Vert w\right\Vert ^{2}-(C_{1}M^{2}+1)\left\vert
w\right\vert ^{2}-C_{3}\left\vert M-N\right\vert ^{2}\left\Vert u\right\Vert
^{2}.
\end{align*}
Hence,
\begin{align*}
&  \frac{d}{dt}\left\vert w\right\vert ^{2}+2\mu\left\Vert w\right\Vert
^{2}+2\left\langle F_{M}(\left\Vert v\right\Vert )\left(
v\text{\textperiodcentered}\nabla\right)  v-F_{N}(\left\Vert u\right\Vert
)\left(  u\text{\textperiodcentered}\nabla\right)  u,w\right\rangle \\
&  =-2\alpha\int_{\Omega}\left(  F_{M}(\left\Vert v\right\Vert ^{\beta
-1})\left\vert v\right\vert ^{\beta-1}v-F_{N}(\left\Vert u\right\Vert
^{\beta-1})\left\vert u\right\vert ^{\beta-1}u\right)  (w)dx+2(f_{v}%
-f_{u},w)\\
&  \leq\frac{\mu}{2}\left\Vert w\right\Vert ^{2}+(C_{1}M^{2}+1)\left\vert
w\right\vert ^{2}+C_{3}|M-N|^{2}||u||^{2}+2\left\vert (f_{v}-f_{u}%
,w)\right\vert .
\end{align*}
We will estimate the term $2\left\langle F_{M}(\left\Vert v\right\Vert
)\left(  v\text{\textperiodcentered}\nabla\right)  v-F_{N}(\left\Vert
u\right\Vert )\left(  u\text{\textperiodcentered}\nabla\right)
u,w\right\rangle $. With an algebraic manipulation we have
\[
\mathcal{LMN}=F_{M}(\left\Vert v\right\Vert )B(v,v)-F_{N}(\left\Vert
u\right\Vert )B(u,u)
\]%
\[
=F_{M}(\left\Vert v\right\Vert )B(w,v)+[F_{M}(\left\Vert v\right\Vert
)-F_{N}(\left\Vert u\right\Vert )]B(u,v)+F_{N}(\left\Vert u\right\Vert
)B(u,w)
\]
and%
\[
2\big |\left\langle F_{M}(\left\Vert v\right\Vert )\left(
v\text{\textperiodcentered}\nabla\right)  v-F_{N}(\left\Vert u\right\Vert
)\left(  u\text{\textperiodcentered}\nabla\right)  u,w\right\rangle \big |
\]%
\[
\leq2|\left\langle F_{M}(\left\Vert v\right\Vert )B(w,v),w\right\rangle
|+2|\left\langle [{F_{M}(\left\Vert v\right\Vert )-F_{N}(\left\Vert
u\right\Vert )]B(u,v),w}\right\rangle |+2|\left\langle F_{N}(\left\Vert
u\right\Vert )B(u,w),w\right\rangle |
\]%
\[
\leq2Mc_{0}||w||^{3/2}|w|^{1/2}+2(\frac{|M-N|}{||u||}+\frac{M||w||}%
{||v||||u||})c_{0}||v||||u|||w|^{1/2}||w||^{1/2}+2Nc_{0}||w||^{3/2}|w|^{1/2}%
\]%
\[
\leq\frac{\mu}{2}||w||^{2}+C_{4}(1+M^{4}+N^{4})|w|^{2}+C_{5}|M-N|^{2}%
||v||^{2},
\]
where we have applied Lemma \ref{IneqFMN} and inequality (\ref{IneqBRomito}).
Also, for the term $(f_{v}-f_{u},v-u)$ we have:%
\[
(f_{v}-f_{u},v-u)\leq C_{6}\left\vert f_{v}-f_{u}\right\vert ^{2}%
+C_{7}\left\vert w\right\vert ^{2}.
\]
We can conclude that
\[
\frac{d}{dt}\left\vert w\right\vert ^{2}+\mu\left\Vert w\right\Vert ^{2}\leq
C_{8}((1+M^{2}+M^{4}+N^{4})\left\vert w\right\vert ^{2}+|M-N|^{2}%
||v||^{2}+|f_{v}-f_{u}|^{2}).
\]
Hence, for $\gamma=\gamma(M,N)=C_{8}(1+N^{2}+M^{4}+N^{4})$ we have
\begin{align*}
|w(t)|^{2}  &  \leq e^{\gamma(t-\tau)}|w_{\tau}|^{2}+C_{8}\left(  \int_{\tau
}^{t}e^{\gamma(t-s)}\left\Vert v(s)\right\Vert ^{2}ds|M-N|^{2}+\int_{\tau}%
^{t}e^{\gamma(t-s)}|f_{v}(s)-f_{u}(s)|^{2}ds\right) \\
&  \leq Ce^{\gamma(t-\tau)}\left(  |w_{\tau}|^{2}+\left(  \int_{\tau}%
^{t}\left\Vert v\right\Vert ^{2}ds\right)  |M-N|^{2}+\int_{\tau}^{t}%
|f_{v}(s)-f_{u}(s)|^{2}ds\right)  .
\end{align*}

\end{proof}

\section{The Kneser property in the weak topology of $H$}

In this section, we prove the main result of this paper, which states that for
$\beta\in\lbrack1,3]$ the attainability set $K_{\tau}(u_{\tau})$ satisfies the
Kneser property, that is, it is compact and connected with respect to the weak
topology of the space $H.$

We start with some auxiliary lemmas. For $\{u_{n}\}\subset C([\tau,T],H_{w})
$, the convergence $u_{n}\rightarrow u$ in $C([\tau,T],H_{w})$ will mean that
for any $v\in H,sup_{t\in\lbrack\tau,T]}|(u_{n}(t)-u(t),v)|\rightarrow0$, as
$n$ $\rightarrow\infty$. We note that $u_{n}\rightarrow u$ in $C([\tau
,T],H_{w})$ if and only if $u\in C([\tau,T],H_{w})$ and for $t_{n}\rightarrow
t_{0}$, $t_{n},t_{0}\in\lbrack\tau,T]$, we have $u_{n}(t_{n})\rightarrow
u(t_{0})$ weakly in $H$.

\begin{lemma}
\label{thm:lemma2} Let $1\leq\beta\leq3$, $f\in L^{2}(\tau,T;H)$ for all
$T>\tau$, $u_{\tau}^{n}\rightarrow u_{\tau}$ strongly in $H$ and let
$u_{n}(\cdot)\in\mathcal{D}_{\tau}(u_{\tau}^{n})$. Then there exists a
subsequence $u_{n_{k}}(\cdot)$ and $u(\cdot)\in\mathcal{D}_{\tau}(u_{\tau})$
such that $u_{n_{k}}\rightarrow u$ in $C([\tau,T],H_{w})$ for all $T>\tau$. If
$\beta=3$, then $u_{n_{k}}\rightarrow u$ in $C([\tau,T],H)$ for all $T>\tau.$
\end{lemma}

\begin{proof}
It follows from the energy inequality \eqref{eq:NSDE-4-Eineq} that each
$u_{n}$ satisfies%
\begin{equation}
|u_{n}(t)|^{2}+\mu\int_{\tau}^{t}||u_{n}(r)||^{2}dr+2\alpha\int_{\tau}%
^{t}||u_{n}(r)||_{\beta+1}^{\beta+1}dr\leq\left\vert u_{\tau}^{n}\right\vert
^{2}+\frac{1}{\mu\lambda_{1}}\int_{\tau}^{t}|f(r)|^{2}dr, \label{eq:NSDE-26}%
\end{equation}
for all $t\in\lbrack\tau,T]$ where $\lambda_{1}$ is the first eigenvalue of
the Stokes operator $A$, so the sequence ${u_{n}}$ is uniformly bounded in
$L^{2}(\tau,T;V)\cap L^{\infty}(\tau,T;H)\cap L^{\beta+1}(\tau,T;(L^{\beta
+1}(\Omega))^{3}$. Hence, there exists a subsequence such that
\begin{align}
u_{n}  &  \rightarrow u\text{ weakly in }L^{2}(\tau,T;V),\label{eq:NSDE-27}\\
u_{n}  &  \rightarrow u\text{ weakly in }L^{\beta+1}(\tau,T;(L^{\beta
+1}(\Omega))^{3},\nonumber\\
u_{n}  &  \rightarrow u\text{ weakly star in }L^{\infty}(\tau,T;H).\nonumber
\end{align}
In addition, in a standard way (see \cite{bib18}) and using $L^{\frac{\beta
+1}{\beta}}(\tau,T,L^{\frac{\beta+1}{\beta}}(\Omega))^{3})\subset L^{4/3}%
(\tau,T,V^{\prime})$ for $1\leq\beta\leq3$, we obtain that
\begin{equation}
\frac{du_{n}}{dt}\rightarrow\frac{du}{dt}\text{ weakly in }L^{\frac{4}{3}%
}(\tau,T;V^{\prime}) \label{eq:NSDE-28}%
\end{equation}
and, by the Compactness Theorem \cite{bib11}, we have
\begin{equation}
u_{n}\rightarrow u\text{ strongly in }L^{2}(\tau,T;H), \label{eq:NSDE-29}%
\end{equation}%
\begin{equation}
u_{n}(t)\rightarrow u(t)\text{ strongly in }H\text{ for a.a. }t\in\left(
\tau,T\right)  . \label{eq:NSDE-30}%
\end{equation}

We can verify that $u$ is a weak solution to (\ref{eq:NSDE-1-intro}) in the
same way as in the proof of the existence of weak solutions by the Galerkin
method (see \cite{bib18}, \cite{bib11}, \cite{bib02}).

The uniform estimate \eqref{eq:NSDE-26} and the compact injection $H\subset
V^{\prime}$ implies that $\{u_{n}(t)\}$ is precompact in $V^{\prime}$ for any
t. Moreover, by (\ref{eq:NSDE-28}) the $u_{n}:[\tau,T]\rightarrow V^{\prime}$
are an equicontinuous family of functions. Hence, by the Ascoli-Arzel\`{a}~
theorem, $u_{n}\rightarrow u$ in $C([\tau,T],V^{\prime})$. By $u\in
C([\tau,T],V^{\prime})\cap L^{\infty}(\tau,T;H)$ we deduce that $u\in
C([\tau,T],H_{w})$ \cite[p.263]{bib18}. Thus, using (\ref{eq:NSDE-26}) and a
standard contradiction argument, we obtain that $u_{n}\rightarrow u$ in
$C([\tau,T],H_{w})$. The fact that this holds for every $T>\tau$ can then be
proved by a diagonal procedure. Combining this convergence with
\eqref{eq:NSDE-27}-\eqref{eq:NSDE-30} and $u_{\tau}^{n}\rightarrow u_{\tau}$
strongly in $H$, we obtain that $u$ satisfies \eqref{eq:NSDE-4-Eineq}. Hence,
$u(\cdot)\in D_{\tau}(u_{\tau})$.

Finally, let $\beta=3$. Using (\ref{eq:NSDE-7}) we obtain that $u_{n},u$
satisfy%
\begin{align*}
\frac{d}{dt}\left\vert u_{n}\right\vert ^{2}  &  \leq\frac{d}{dt}\left\vert
u_{n}\right\vert ^{2}+\mu||u_{n}||^{2}+2\alpha\left\Vert u_{n}\right\Vert
_{\beta+1}^{\beta+1}\leq\frac{|f|^{2}}{\mu\lambda_{1}},\\
\frac{d}{dt}\left\vert u\right\vert ^{2}  &  \leq\frac{d}{dt}\left\vert
u\right\vert ^{2}+\mu||u||^{2}+2\alpha\left\Vert u\right\Vert _{\beta
+1}^{\beta+1}\leq\frac{|f|^{2}}{\mu\lambda_{1}},
\end{align*}
so
\begin{align*}
|u_{n}(t)|^{2}  &  \leq\left\vert u_{n}(s)\right\vert ^{2}+\frac{1}{\mu
\lambda_{1}}\int_{s}^{t}|f(r)|^{2}dr,\\
|u(t)|^{2}  &  \leq\left\vert u(s)\right\vert ^{2}+\frac{1}{\mu\lambda_{1}%
}\int_{s}^{t}|f(r)|^{2}dr,\ \forall t\geq s\geq\tau.
\end{align*}
We define the maps $J_{n}(t)=|u_{n}(t)|^{2}-\frac{1}{\mu\lambda_{1}}\int%
_{\tau}^{t}|f(r)|^{2}dr.$ By Corollary \ref{DerivCor} they are continuous and
by the last inequalities they are non-increasing.

Let $t_{n}\rightarrow t_{0}>\tau$. It follows from (\ref{eq:NSDE-30}) that
$J_{n}(t)\rightarrow J(t)$ for a.a. $t$. As $u\in C([\tau,T],H)$, for any
$\delta>0$ there are $t_{k}<t_{0}$ and $N(t_{k})$ such that $\left\vert
J(t_{k})-J(t_{0})\right\vert <\delta/2$ and $\left\vert J_{n}(t_{k}%
)-J(t_{k})\right\vert <\delta/2$ for $n\geq N.$ For $n$ great enough we have
$t_{n}>t_{k}$, so by $J_{n}(t_{n})\leq J_{n}(t_{k})$ we derive that%
\[
J_{n}(t_{n})-J(t_{0})\leq\left\vert J_{n}(t_{k})-J(t_{k})\right\vert
+\left\vert J(t_{k})-J(t_{0})\right\vert <\delta\text{.}%
\]
Hence, $\lim\sup\ |u_{n}(t_{n})|\leq|u(t_{0})|$. Combining this with
$|u(t_{0})|\leq\lim\inf\ |u_{n}(t_{n})|$ (by $u_{n}(t_{n})\rightarrow
u(t_{0})$ weakly in $H$), we get $|u_{n}(t_{n})|\rightarrow|u(t_{0})|$ and
then $u_{n}(t_{n})\rightarrow u(t_{0})$ strongly in $H$.

Let $t_{n}\rightarrow\tau$. Repeating the above argument with $t_{k}%
=\tau=t_{0}$, we obtain that $u_{n}(t_{n})\rightarrow u(\tau)$ strongly in $H$.
\end{proof}

\begin{lemma}
\label{thm:lemma3} Let $1\leq\beta\leq3$ and $u_{\tau}^{n}\rightarrow u_{\tau
}$ weakly in $H$ and $f^{n}\rightarrow f$ weakly in $L^{2}(\tau,T;H)$ for all
$T>\tau.$ If $u_{n}(\cdot)$ is the unique weak solution to
\eqref{eq:NSDE-8-GMNSED} with $u_{n}(\tau)=u_{\tau}^{n}$, then, $u_{n}%
\rightarrow u$ in $C([\tau,T],H_{w})$ for all $T\geq\tau$ where $u(\cdot)$ is
the unique weak solution to \eqref{eq:NSDE-8-GMNSED} with $u(\tau)=u_{\tau}$.
Moreover, $u_{n}\rightarrow u$ in $C([t_{0},T,H)$ for all $T>t_{0}>\tau$. If
$u_{\tau}^{n}\rightarrow u_{\tau}$ strongly, then $u_{n}\rightarrow u$ in
$C([\tau,T];H)$ for all $T>\tau$.
\end{lemma}

\begin{proof}
It follows from \eqref{eq:NSDE-N-Eineq} that
\begin{equation}
|u_{n}(t)|^{2}+\mu\int_{\tau}^{t}||u_{n}(r)||^{2}dr+2\alpha\int_{\tau}%
^{t}F_{N}(||u_{n}(t)||^{\beta-1})|u_{n}(r)|_{\beta+1}^{\beta+1}dr\leq|u_{\tau
}^{n}|^{2}+\frac{1}{\mu\lambda_{1}}\int_{\tau}^{t}|f^{n}(r)|^{2}dr
\label{NSDE-32F}%
\end{equation}
for all $t>\tau$. Thus, \eqref{NSDE-32F} implies that $u_{n}$ is bounded in
$L^{2}(\tau,T;V)\cap L^{\infty}(\tau,T;H).$ In addition, using
\eqref{eq:NSDE-9} and the fact that $F_{N}(||u_{n}||^{\beta-1})|u_{n}%
|^{\beta-1}u$ is bounded in $L^{\frac{\beta+1}{\beta}}(\tau,T;(L^{\frac
{\beta+1}{\beta}}(\Omega))^{3})\subset L^{4/3}(\tau,T;V^{\prime})$, we see
that the $\frac{d}{dt}u_{n}(\cdot)$ are bounded in $L^{\frac{4}{3}}%
(\tau,T,V^{\prime})$. The existence of a subsequence converging to the weak
solution to \eqref{eq:NSDE-8-GMNSED} with $u(\tau)=u_{\tau}$ is proved in the
same way as in Theorems \ref{thm:Th03-primero}-\ref{thm:Th04-segundo}. The
convergences in $C([\tau,T],H_{w})$, $C([t_{0},T,H)$ and $C([\tau,T,H)$ are
proved in the same way as in Lemma \ref{thm:lemma2}.

As the limit $u$ is the same for every subsequence, a standard contradiction
argument shows that the whole sequence converges.
\end{proof}

\bigskip

Now, we are ready to establish the Kneser property of the attainability set
$K_{\tau}(u_{\tau}).$

\begin{theorem}
\label{thm:Th5} Let $1\leq\beta\leq3$, $f\in L^{2}(\tau,T;H)$ for all $T>\tau
$. Then, for all $t\geq\tau$ and $u_{\tau}\in H$ the attainability set
$K_{\tau}(u_{\tau})$ is compact and connected with respect to the weak
topology on $H$.
\end{theorem}

\begin{proof}
The proof that $K_{\tau}(t)$ is weakly compact for any $t\geq\tau$ is an
immediate consequence of Lemma \ref{thm:lemma2}. We will prove that it is
connected with respect to the weak topology of $H.$

The case $t=\tau$ is obvious, so suppose that the set $K_{\tau^{\ast}}%
(u_{\tau})$ is not weakly connected for some $t^{\ast}>\tau$. Then there exist
two weakly compact sets $A_{1}$, $A_{2}$ of $H$ with $A_{1}\cap A_{2}%
=\emptyset$, such that $A_{1}\cup A_{2}=K_{\tau^{\ast}}(u_{\tau})$. Let
$u_{1}(\cdot),u_{2}(\cdot)\in D_{\tau}(u_{\tau})$ be such that $u_{1}(t^{\ast
})\in U_{1}$ and $u_{2}(t^{\ast})\in U_{2}$, where $U_{1}$ and $U_{2}$ are
disjoint weakly open neighbourhoods of $A_{1}$ and $A_{2}$, respectively.

Now, for $i=1,2$ and $\gamma\geq\tau$, let $U_{i}^{N}(t,\gamma)$ be equal to
$\{u_{i}(t)\}$ if $t\in\lbrack\tau,\gamma]$ and if $t\geq\gamma$ let
$U^{N}(t,\gamma)$ be the solution to the problem%
\begin{equation}%
\begin{cases}
u_{t}-\mu\Delta u+F_{N}(\left\Vert u\right\Vert )(u\cdot\nabla)u+\nabla
p+\alpha F_{N}(\left\Vert u\right\Vert ^{\beta-1})|u|^{\beta-1}u=f,\\
\operatorname{div}\text{ }u=0,\\
u|_{\partial\Omega}=0,\\
u(\gamma)=u_{i}(\gamma).
\end{cases}
\label{eq:NSDE-33}%
\end{equation}

By Theorem \ref{thm:Th04-segundo} we know that $U_{i}^{N}(t,\gamma)$ is a well
defined function for any $N>0$ and $t,\gamma\geq\tau.$ Now, we shall prove
that the maps $\gamma\rightarrow U_{i}^{N}(t,\gamma)$ are continuous in
$H_{w}$. Let $\gamma\rightarrow\gamma_{0}$.

Consider first the case where $\gamma>\gamma_{0}$, with $\gamma\rightarrow
\gamma_{0}$. If $t\leq\gamma_{0}<\gamma$, then $U_{i}^{N}(t,\gamma
)=\{u_{i}(t)\}=U_{i}^{N}(t,\gamma_{0})$. If $t>\gamma_{0}$, we can assume that
$t>\gamma$, so that $U_{i}^{N}(t,\gamma)$ is the weak solution to problem
\eqref{eq:NSDE-33} with $U_{i}^{N}(\gamma,\gamma)=u_{i}(\gamma)$. Since
$u_{i}\in C([\tau,\infty],H_{w})$, we have $u_{i}(\gamma)\rightarrow
u_{i}(\gamma_{0})$ weakly in $H$. We will omit the index $i$. If
$\gamma\rightarrow U^{N}(t,\gamma)$ were not continuous for $\gamma_{0}$ when
$\gamma\rightarrow\gamma_{0}$ then there would exist a neighbourhood $O$ of
$U^{N}(t,\gamma_{0})$ in the weak topology and sequences $\gamma_{j}%
>\gamma_{0}$ with $\gamma_{j}\rightarrow\gamma_{0}$ and $U^{N}(t,\gamma
_{j})\notin O$ for all $j$. Denote $\bar{y}(t)=U^{N}(t+\gamma_{j}-\gamma
_{0},\gamma_{j})$ for $t\geq\gamma_{0}$. Then $\bar{y}_{j}(t)$ is the weak
solution to \eqref{eq:NSDE-33} with $\bar{y}_{j}(\gamma_{0})=u(\gamma_{j})$,
but with $f(\cdot)$ replaced by $f(\cdot+\gamma_{j}-\gamma_{0})$. Since
$f(\cdot+\gamma_{j}-\gamma_{0})\rightarrow f(\cdot)$ strongly in $L^{2}%
(\gamma_{0},T;H)$ for all $T>\gamma_{0}$ \cite{Gaevski}, it follows from Lemma
\ref{thm:lemma3} that $\bar{y}_{j}(\cdot)\rightarrow\bar{y}$ in $C([\gamma
_{0},T],H_{w})$ for all $T>\gamma_{0}$, where $\bar{y}(\cdot)$ is the weak
solution to \eqref{eq:NSDE-33} with $\bar{y}(\gamma_{0})=u(\gamma_{0})$. But
then $U^{N}(t,\gamma_{j})=\bar{y}_{j}(t-\gamma_{j}+\gamma_{0})\rightarrow
\bar{y}(t)=U^{N}(t,\gamma_{0})$ in $H_{w}$ and we have obtained a
contradiction. Consider further the case where $\gamma<\gamma_{0}$ with
$\gamma\rightarrow\gamma_{0}$. If $t<\gamma_{0}$, then can assume that
$t<\gamma$ and so $U^{N}(t,\gamma)=\{u(t)\}=U^{N}(t,\gamma_{0})$. If
$t\geq\gamma_{0}>\gamma$, then we essentially repeat the proof above.

For any fixed $T>t^{\ast}$, we write:
\[
\gamma(\rho)=\left\{
\begin{array}
[c]{c}%
\tau-\left(  T-\tau\right)  \rho\text{ if }\rho\in\lbrack-1,0],\\
\tau+\left(  T-\tau\right)  \rho\text{ if }\rho\in\lbrack0,1],
\end{array}
\right.
\]
and define the family of functions $\Phi^{N}(\rho):[t,T]\rightarrow H$ with
$\rho\in\lbrack-1,1]$ by%
\[
\Phi^{N}(\rho)(t)=\left\{
\begin{array}
[c]{c}%
U_{1}^{N}(t,\gamma(\rho))\text{ if }\rho\in\lbrack-1,0],\\
U_{2}^{N}(t,\gamma(\rho))\text{ if }\rho\in\lbrack0,1].
\end{array}
\right.
\]

Then $\Phi^{N}(-1)(t)=U_{1}^{N}(t,T)=u_{1}(t)$ and $\Phi^{N}(1)(t)=U_{2}%
^{N}(t,T)=u_{2}(t)$. Moreover the mapping $\rho\rightarrow\Phi^{N}(\rho)(t)$
is continuous for any fixed $N\geq1$ and $t\in\lbrack\tau,T]$ (because
$U_{1}^{N}(t,\tau)=U_{2}^{N}(t,\tau)$). $U_{\rho\in\lbrack-1,1]}\Phi^{N}%
(\rho)(t)$ is connected in $H_{w}$ for each fixed $N\geq1$ and $t\in
\lbrack\tau,T]$. In particular, $U_{\rho\in\lbrack-1,1])}\Phi^{N}%
(\rho)(t^{\ast})$ is connected. Indeed, we define $F:[-1,1]\rightarrow H$ by
$F(\rho):=\Phi^{N}(\rho)(t)$ and $A=[-1,1]$. If $F(A)$ were not connected,
then there would exist open sets $O_{1}$ and $O_{2}$ in $H_{w}$ with
$O_{1}\cap O_{2}=\emptyset$ such that $F(A)\cap O_{i}\neq\emptyset$, $i=1,2$,
and $F(A)\subset O_{1}\cup O_{2}$. Denote $M_{i}=\{\rho\in A:F(\rho)\in
O_{i}\}$. It is clear that $M_{1}\cup M_{2}=A$. Moreover, $M_{1}\cap
M_{2}=\emptyset$ and $M_{i}\neq\emptyset$ for $i=1,2$. Since $F$ is continuous
in $H_{w}$, it is easy to see that $V_{i}=\{\rho\in A:F(\rho)\subset O_{i}\}$
are open sets in $H_{w}$ for $i=1,2$. Then $M_{i}\subset V_{i}$ and $V_{1}\cap
V_{2}=\emptyset$, which contradicts the fact that $A$ is a connected set.

Since $\Phi^{N}(-1)(t^{\ast})=u_{1}(t^{\ast})\in U_{1}$ and $\Phi
^{N}(1)(t^{\ast})={u_{2}(t^{\ast})\in U_{2}}$, there exists a $\rho_{N}%
\in(-1,1)$ such that $\Phi^{N}(\rho_{N})(t^{\ast})\notin U_{1}\cup U_{2}$.
There are infinitely many $N$ such that the $\Phi^{N}(\rho_{N})(t^{\ast})$ are
equal either to $U_{1}^{N}(t^{\ast},\gamma(\rho_{N}))$ or $U_{2}^{N}(t^{\ast
},\gamma(\rho_{N}))$. Suppose, for example, that this is true for $U_{1}%
^{N}(t^{\ast},\gamma(\rho_{N}))$.

We put $u^{N}(t)=U_{1}^{N}(t,\gamma(\rho_{N}))$ and
\[
\bar{F}_{N}(||u^{N}(r)||^{\beta-1})=\left\{
\begin{array}
[c]{c}%
1\text{ if }t\in\lbrack\tau,\gamma(\rho_{N})],\\
{F}_{N}(||u^{N}(r)||^{\beta-1})\text{ if }t\in(\gamma(\rho_{N}),T].
\end{array}
\right.
\]
The function $u_{1}(t)$ satisfies the energy inequality
\eqref{eq:NSDE-4-Eineq} and the weak solutions to system \eqref{eq:NSDE-33}
satisfy the energy inequality (\ref{eq:NSDE-N-Eineq}). We will see that
$u^{N}(t)$ also satisfies the energy inequality (\ref{eq:NSDE-N-Eineq}) but
replacing $F_{N}$ by $\bar{F}_{N}$. The case $t\leq\gamma(\rho_{N})$ is
trivial, so let us consider the case $t\geq s\geq\gamma(\rho_{N})$. Then
\[
\frac{1}{2}|u^{N}(t)|^{2}+\mu\int_{\tau}^{t}||u^{N}(r)||^{2}dr+\alpha
\int_{\tau}^{t}\bar{F}_{N}(||u^{N}(r)||^{\beta-1})|u^{N}(r)|_{\beta+1}%
^{\beta+1}dr-\int_{\tau}^{t}(f(r),u^{N}(r))dr
\]%
\[
=\frac{1}{2}|u_{N}(t)|^{2}+\mu\int_{\gamma(\rho_{N})}^{t}||u^{N}%
(r)||^{2}dr+\alpha\int_{\gamma(\rho_{N})}^{t}F_{N}(||u^{N}(r)||^{\beta
-1})|u^{N}(r)|_{\beta+1}^{\beta+1}dr-\int_{\gamma(\rho_{N})}^{t}%
(f(r),u^{N}(r))dr
\]%
\[
+\mu\int_{\tau}^{\gamma(\rho_{N})}||u_{1}(r)||^{2}dr+\alpha\int_{\tau}%
^{\gamma(\rho_{N})}|u_{1}(r)|_{\beta+1}^{\beta+1}dr-\int_{\tau}^{\gamma
(\rho_{N})}(f(r),u_{1}(r))dr
\]%
\[
\leq\frac{1}{2}|u_{N}(s)|^{2}+\mu\int_{\gamma(\rho_{N})}^{s}||u^{N}%
(r)||^{2}dr+\alpha\int_{\gamma(\rho_{N})}^{s}F_{N}(||u^{N}(r)||^{\beta
-1})|u^{N}(r)|_{\beta+1}^{\beta+1}dr-\int_{\gamma(\rho_{N})}^{s}%
(f(r),u^{N}(r))dr
\]%
\[
+\mu\int_{\tau}^{\gamma(\rho_{N})}||u_{1}(t)||^{2}dr+\alpha\int_{\tau}%
^{\gamma(\rho_{N})}|u_{1}(r)|_{\beta+1}^{\beta+1}dr-\int_{\tau}^{\gamma
(\rho_{N})}(f(r),u_{1}(r))dr
\]%
\[
=\frac{1}{2}|u^{N}(s)|^{2}+\mu\int_{\tau}^{s}||u^{N}(r)||^{2}dr+\alpha
\int_{\tau}^{s}\bar{F}_{N}(||u^{N}(r)||^{\beta-1})|u^{N}(r)|_{\beta+1}%
^{\beta+1}dr-\int_{\tau}^{s}(f(r),u^{N}(r))dr.
\]
Finally, for $t\geq\gamma(\rho_{N})>s\geq\tau$ , we have
\[
\frac{1}{2}|u^{N}(t)|^{2}+\mu\int_{\tau}^{t}||u^{N}(r)||^{2}dr+\alpha
\int_{\tau}^{t}\bar{F}_{N}(||u^{N}(r)||^{\beta-1})|u^{N}(r)|_{\beta+1}%
^{\beta+1}dr-\int_{\tau}^{t}(f(r),u^{N}(r))dr
\]%
\[
=\frac{1}{2}|u_{N}(t)|^{2}+\mu\int_{\gamma(\rho_{N})}^{t}||u^{N}%
(r)||^{2}dr+\alpha\int_{\gamma(\rho_{N})}^{t}F_{N}(||u^{N}(r)||^{\beta
-1}|u^{N}(r)|_{\beta+1}^{\beta+1}dr-\int_{\gamma(\rho_{N})}^{t}(f(r),u^{N}%
(r))dr
\]%
\[
+\mu\int_{\tau}^{\gamma(\rho_{N})}||u_{1}(r)||^{2}dr+\alpha\int_{\tau}%
^{\gamma(\rho_{N})}|u_{1}(r)|_{\beta+1}^{\beta+1}dr-\int_{\tau}^{\gamma
(\rho_{N})}(f(r),u_{1}(r))dr
\]%
\[
\leq\frac{1}{2}|u_{1}(\gamma(\rho_{N}))|^{2}+\mu\int_{\tau}^{\gamma(\rho_{N}%
)}||u_{1}(r)||^{2}dr+\alpha\int_{\tau}^{\gamma(\rho_{N})}|u_{1}(r)|_{\beta
+1}^{\beta+1}dr-\int_{\tau}^{\gamma(\rho_{N})}(f(r),u_{1}(r))dr
\]%
\[
\leq\frac{1}{2}|u_{1}(s)|^{2}+\mu\int_{\tau}^{s}||u_{1}(r)||^{2}dr+\alpha
\int_{\tau}^{s}|u_{1}(r)|_{\beta+1}^{\beta+1}dr-\int_{\tau}^{s}(f(r),u_{1}%
(r))dr.
\]
It then follows that
\begin{equation}
|u^{N}(t)|^{2}+\mu\int_{\tau}^{t}||u^{N}(r)||dr+2\alpha\int_{\tau}^{t}\bar
{F}_{N}(||u^{N}(r)||^{\beta-1})|u^{N}(r)|_{\beta+1}^{\beta+1}dr\leq|u_{\tau
}|^{2}+\frac{1}{\mu\lambda_{1}}\int_{\tau}^{t}|f(r)|^{2}dr \label{Inequn}%
\end{equation}
so that $u^{N}$ is a bounded sequence in $L^{\infty}(\tau,T;H)\cap L^{2}%
(\tau,T;V)$. In addition, by a standard inequality (see Proposition 9.2 in
\cite{bib14}) and the inequality \eqref{eq:NSDE-9} for $B_{N}$, we see that
the mapping
\[
\bar{B}_{N}(u^{N}(t),u^{N}(t))=\left\{
\begin{array}
[c]{c}%
B(u^{N}(t),u^{N}(t))\text{ if }t\in\lbrack\tau,\gamma(\rho_{N})],\\
B_{N}(u^{N}(t),u^{N}(t))\text{ if }t\in(\gamma(\rho_{N}),T].
\end{array}
\right.
\]
is bounded in the space $L^{4/3}(t,T;V^{\prime})$. Since $u^{N}$ satisfies the
evolution equation
\[
\frac{d}{dt}u^{N}+\mu Au^{N}+\bar{B}_{N}(u^{N}(t),u^{N}(t))+\alpha\bar{F}%
_{N}(||u^{N}||^{\beta-1})|u^{N}|^{\beta-1}u^{N}=f,
\]
it follows that the $\frac{d}{dt}u^{N}$ are bounded in $L^{4/3}(\tau
,T;V^{\prime})$ because
\[
\int_{\tau}^{t}\int_{\Omega}\Big|\bar{F}_{N}(||u^{N}||^{\beta-1}%
)|u^{N}|^{\beta-1}{u^{N}}\Big|^{\frac{\beta+1}{\beta}}dxdt\leq\int_{\tau}%
^{t}\int_{\Omega}\bar{F}_{N}(||u^{N}||^{\beta-1})|u^{N}|^{\beta+1}dxdt
\]%
\[
\leq\int_{\tau}^{t}\bar{F}_{N}(||u^{N}||^{\beta-1})|u^{N}|_{\beta+1}^{\beta
+1}dt\leq C<\infty.
\]
Then, by the Compactness Theorem \cite{bib11}, there exists a subsequence
$N_{j}\rightarrow\infty$ and a $u\in L^{\infty}(\tau,T;H)\cap L^{2}%
(\tau,T;V)\cap L^{\beta+1}(\tau,T;(L^{\beta+1}(\Omega))^{3})$ such that
\[
u^{N_{j}}\rightharpoonup u\text{ weakly in }L^{2}(\tau,T;V)\text{,}%
\]%
\[
u^{N_{j}}\rightharpoonup u\text{ weakly-star in }L^{\infty}(\tau,T;V)\text{,}%
\]%
\[
u^{N_{j}}\rightarrow u\text{ strongly in }L^{2}(\tau,T;H)\text{,}%
\]%
\[
\frac{du^{N_{j}}}{dt}\rightharpoonup\frac{du}{dt}\text{ weakly in }L^{\frac
{4}{3}}(\tau,T;V)^{\prime}\text{.}%
\]

We have to show that $u$ is a weak solution to (\ref{eq:NSDE-1-intro}) on
$[\tau,T]$. For this, it is enough to check that
\[
(u(t),w)+\mu\int_{\tau}^{T}((u(s),w))ds+\int_{\tau}^{T}b(u(s),u(s),w)ds+\alpha
\int_{\tau}^{T}(|u(s)|^{\beta-1}u,w)=(u_{\tau},w)+\int_{\tau}^{T}\langle
f(s),w\rangle ds
\]
for any $w\in V$.

This is done in the same way as for the convergence to a weak solution of
Galerkin approximations for the Navier-Stokes system. Differences arise only
with the nonlinear terms $\bar{B}_{N_{j}}$ and $\overline{G}_{N_{j}},$ where
\[
\overline{G}_{N}(u^{N}(t),u^{N}(t))=\left\{
\begin{array}
[c]{c}%
G(u^{N}(t),u^{N}(t))\text{ if }t\in\lbrack\tau,\gamma(\rho_{N})],\\
G_{N}(u^{N}(t),u^{N}(t))\text{ if }t\in(\gamma(\rho_{N}),T].
\end{array}
\right.
\]

For the term $\overline{G}_{N_{j}}$ it is enough to prove that
\begin{equation}
\bar{F}_{N_{j}}(||u^{N_{j}}||^{\beta-1})|u^{N_{j}}|^{\beta-1}u^{N_{j}%
}\rightarrow|u|^{\beta-1}u\text{ weakly in }L^{q}(\tau,T;\left(  L^{q}%
(\Omega)\right)  ^{3})\text{,} \label{ConvergGNA}%
\end{equation}%
\begin{equation}
\bar{F}_{N_{j}}(||u^{N_{j}}||^{\beta-1})u^{N_{j}}\rightarrow u\text{ weakly in
}L^{\beta+1}(\tau,T;\left(  L^{\beta+1}(\Omega)\right)  ^{3},
\label{ConvergGNb}%
\end{equation}
where $q=(\beta+1)^{\prime}=\frac{\beta+1}{\beta}$. We know from
(\ref{Inequn}) that
\[
\int_{\tau}^{T}\bar{F}_{N_{j}}(||u^{N_{j}}(s)||^{\beta-1})|u^{N_{j}%
}(s)|_{\beta+1}^{\beta+1}ds\leq C,
\]
so using $0\leq\bar{F}_{N_{j}}\left(  r\right)  \leq1$ we get
\[
\int_{\tau}^{T}|\bar{F}_{N_{j}}(||u^{N_{j}}(s)||^{\beta-1})u^{N_{j}%
}(s)|_{\beta+1}^{\beta+1}ds\leq\int_{\tau}^{T}\bar{F}_{N_{j}}(||u^{N_{j}%
}(s)||^{\beta-1})|u^{N_{j}}(s)|_{\beta+1}^{\beta+1}ds\leq C,
\]%
\[
\int_{\tau}^{T}|\bar{F}_{N_{j}}(||u^{N_{j}}(s)||^{\beta-1})|u^{N_{j}%
}(s)|^{\beta-1}u^{N_{j}}(s)|_{q}^{q}\leq\int_{\tau}^{T}\bar{F}_{N_{j}%
}(||u^{N_{j}}(s)||^{\beta-1})|u^{N_{j}}(s)|_{\beta+1}^{\beta+1}ds\leq C,
\]
and then%
\begin{equation}
\bar{F}_{N_{j}}(||u^{N_{j}}||^{\beta-1})u^{N_{j}}\text{ is bounded in
}L^{\beta+1}(\tau,T;\left(  L^{\beta+1}(\Omega)\right)  ^{3}),
\label{FNumBoundedA}%
\end{equation}%
\begin{equation}
\bar{F}_{N_{j}}(||u^{N_{j}}||^{\beta-1})|u^{N_{j}}|^{\beta-1}u^{N_{j}}\text{
is bounded in }L^{q}(\tau,T;\left(  L^{q}(\Omega)\right)  ^{3}),
\label{GmBoundedLqA}%
\end{equation}%
\[
\bar{F}_{N_{j}}(||u^{N_{j}}||^{\beta-1})u^{N_{j}}\rightarrow\chi_{1}\text{
weakly in }L^{\beta+1}(\tau,T;\left(  L^{\beta+1}(\Omega)\right)  ^{3}),
\]%
\[
\bar{F}_{N_{j}}(||u^{N_{j}}||^{\beta-1})|u^{N_{j}}|^{\beta-1}u^{N_{j}%
}\rightarrow\chi_{2}\text{ weakly in }L^{q}(\tau,T;\left(  L^{q}%
(\Omega)\right)  ^{3}),
\]
for some $\chi_{1},\chi_{2}$. We need to show that%
\begin{equation}
\chi_{1}=u,\ \chi_{2}=|u|^{\beta-1}u. \label{X}%
\end{equation}
Since $u^{N_{j}}\rightarrow u$ a.e. in $(\tau,T)\times\Omega$ and $F_{N_{j}%
}(||u^{N_{j}}(s)||)\rightarrow\bar{F}_{N_{j}}(||u(s)||)$ for a.a. $s\in\left(
\tau,T\right)  $, we deduce that
\begin{align*}
\bar{F}_{N_{j}}(||u^{N_{j}}||^{\beta-1})||u^{N_{j}}|^{\beta-1}u^{N_{j}}  &
\rightarrow|u|^{\beta-1}u,\\
\bar{F}_{N_{j}}(||u^{N_{j}}||^{\beta-1})u^{N_{j}}  &  \rightarrow u\text{ a.e.
in }(\tau,T)\times\Omega.
\end{align*}
By (\ref{FNumBoundedA}), (\ref{GmBoundedLqA}), equality (\ref{X}) follows from
Lemma 1.3 in \cite{bib11}. We have proved (\ref{ConvergGNA})-(\ref{ConvergGNb}).

Further, we have to show that
\begin{equation}
\int_{\tau}^{t}\langle\bar{B}_{N_{j}}(u^{N_{j}}(r),u^{N_{j}}(r)),w\rangle
dr\rightarrow\int_{\tau}^{t}b(u(r),u(r),w)dr \label{eq:NSDE-38}%
\end{equation}
for all $t\in\lbrack\tau,T]$ and $w\in D(A)$. Note that
\[
\bar{B}_{N_{j}}(u^{N_{j}}(r),u^{N_{j}}(r))=\bar{F}_{N_{j}}(||u^{N_{j}%
}(t)||)B(u^{N_{j}}(t),u^{N_{j}}(t)).
\]

\begin{lemma}
\label{thm:lemma4} $\bar{F}_{N_{j}}(||u^{N_{j}}(t)||)\rightarrow1$ in
$L^{p}(\tau,T;\mathbb{R})$ when $N_{j}\rightarrow\infty$ for all $p>1$.
\end{lemma}

\begin{proof}
It is proved in Lemma 12 from \cite{bib28} that $F_{N_{j}}(||u^{N_{j}%
}(t)||\rightarrow1$ in $L^{p}(\tau,T;\mathbb{R})$ as $N_{j}\rightarrow\infty$
for all $p>1$. The result then follows from
\[
\int_{\tau}^{T}|\bar{F}_{N_{j}}(||u^{N_{j}}(t)||)-1|^{p}dt\leq\int%
_{\gamma(\rho_{N_{j}})}^{T}|F_{N_{j}}(||u^{N_{j}}(t)||)-1|^{p}dt\rightarrow0.
\]

\end{proof}

\bigskip

Arguing as in the proof of the existence of weak solutions for the
Navier-Stokes System \cite{bib11}, it can be shown that
\[
\int_{\tau}^{T}b(u^{N_{j}}(r),u^{N_{j}}(r),w)dr\rightarrow{\int_{\tau}%
^{T}b(u(r),u(r),w)dr}%
\]
Then
\[
\int_{\tau}^{T}[\bar{F}_{N_{j}}(||u^{N_{j}}(t)||)b(u^{N_{j}}(r),u^{N_{j}%
}(r),u^{N_{j}}(r),w)dr-b(u(r),u(r),w)]dr
\]%
\[
=\int_{\tau}^{T}[\bar{F}_{N_{j}}(||u^{N_{j}}(t)||)b(u^{N_{j}}(r),u^{N_{j}%
}(r),u^{N_{j}}(r),w)dr-b(u^{N_{j}}(r),u^{N_{j}}(r),w)]dr
\]%
\[
+\int_{\tau}^{T}[b(u^{N_{j}}(r),u^{N_{j}}(r)w)-b(u(r),u(r),w)]dr
\]
and the second integral converges to 0. For the first integral, we use the
inequality
\[
\big |\int_{\tau}^{T}(\bar{F}_{N_{j}}(||u^{N_{j}}(t)||)-1)b(u^{N_{j}%
}(r),u^{N_{j}}(r),w)dr\big|^{2}%
\]%
\[
\leq\int_{\tau}^{T}(\bar{F}_{N_{j}}(|u^{N_{j}}(t)||)-1)^{2}dr\int_{\tau}%
^{T}\big |b(u^{N_{j}}(r),u^{N_{j}}(r),w)\big |^{2}dr
\]
and the estimate
\[
\int_{\tau}^{T}\big |b(u^{N_{j}}(r),u^{N_{j}}(r),w)\big |^{2}dr\leq
C||w||\cdot|Aw|\int_{\tau}^{T}||u^{N_{j}}(r)||^{2}|u^{N_{j}}(r)|^{2}\leq C
\]
which follows from (9.2) in \cite{bib14}.

Hence, the first integral above also converges to $0$, which gives the limit
\eqref{eq:NSDE-38} for all $w\in D(A)$. Finally, a density argument gives the
limit for all $w\in V$.

It follows that $u(\cdot)$ is a weak solution to (\ref{eq:NSDE-1-intro})
defined on $[\tau,T]$. In fact, the function $u^{N}(t)$ can be defined also in
$[\tau,2T],\ [\tau,3T],$ etc., and by a standard diagonal argument we obtain
that $u(\cdot)$ is a globally defined weak solution. In addition, since
$u^{N}(\cdot)$ satisfies the energy inequality \eqref{eq:NSDE-N-Eineq}, it is
easy to show that $u(\cdot)$ satisfies \eqref{eq:NSDE-4-Eineq}. Hence,
$u(\cdot)\in D_{\tau}(u_{\tau})$. Finally, we have to show that
\begin{equation}
u^{N_{j}}(t^{\ast})\rightarrow u(t^{\ast})\text{ weakly in }H.
\label{eq:ultima}%
\end{equation}
For any $t\in\lbrack\tau,T]$, the uniform estimate $|u^{N_{j}}(t)|\leq C$ and
the compact embedding $H\subset V^{\prime}$ implies that the sequence
$\{u^{N_{j}}(t)\}$ is precompact in $V^{\prime}$. In addition, the sequence
$\frac{du^{N_{j}}}{dt}$ is bounded in $L^{4/3}([0,T];V^{\prime})$. Hence by
Ascoli-Arzel\`{a} theorem, the sequence $\{u^{N_{j}}\}$ is relatively compact
in $C([\tau,T],V^{\prime})$. Then, passing to a subsequence, we have
$u^{N_{j}(t^{\ast})}\rightarrow u(t^{\ast})$ strongly in $V^{\prime}$. From
this and the estimate $|u^{N_{j}}(t)|\leq C$, we obtain \eqref{eq:ultima} by a
standard argument.

Hence, $u(t^{\ast})\notin U_{1}\cup U_{2}$ and $u(t^{\ast})\in K_{t^{\ast}%
}(u_{\tau})$, which is a contradiction. Thus, $K_{t^{\ast}}(u_{\tau})$ must be
connected with respect to the weak topology of $H$, which completes the proof
of Theorem \ref{thm:Th5}.
\end{proof}

\section{The Kneser property in the strong topology of $H$ for $\beta
=3,\ 4\mu\alpha<1$}

In the particular case where $\beta=3,\ 4\mu\alpha<1$, we are able to prove
the Kneser property of the attainability set with respect to the strong
topology of the space $H.$

\begin{theorem}
\label{Thm:ThCON} Let $\beta=3$, $4\mu\alpha<1$ and $f\in L^{2}(\tau,T;H)$ for
all $T>\tau$. Then, for each $t\geq\tau$ and $u_{\tau}$ in $H$, the
attainability set $K_{t}(u_{\tau})$ is compact and connected with respect to
the strong topology of $H$.
\end{theorem}

\begin{proof}
The proof that $K_{t}(u_{\tau})$ is compact for any $t\geq\tau$ is an
immediate consequence of Lemma \ref{thm:lemma2}. We will prove that it is
connected with respect to the strong topology of $H.$

The case $t=\tau$ is obvious. Suppose that for some $t^{\ast}>\tau$ the set
$K_{t^{\ast}}(u_{\tau})$ is not connected. Then there exists two non-empty
compact sets $A_{1}$, $A_{2}$, such that $A_{1}\cup A_{2}=K_{t^{\ast}}%
(u_{\tau})$ and $A_{1}\cap A_{2}=\emptyset$. Let $u_{1}(\cdot),u_{2}(\cdot)\in
D_{t}(u_{\tau})$ be such that $u_{1}(t^{\ast})\in U_{1}$ and $u_{2}(t^{\ast
})\in U_{2}$, where $U_{1}$ and $U_{2}$ are disjoint open neighbourhoods of
$A_{1}$ and $A_{2}$, respectively.

We are going to use the same definitions as in Theorem \ref{thm:Th5}. Let
$u_{1}(t)$ and $u_{2}(t)$ be weak solutions of (\ref{eq:NSDE-3-weaksol}). The
functions $U_{i}^{N}(t,\gamma),\gamma(\rho),\phi^{N}(\rho)(t)$ are the same as
above and $U_{i}^{N}(t,\gamma)$ are continuous on $\gamma$ for each $N\geq0$
and $t\in\lbrack\tau,T]$. The proof is essentially the same as in Theorem
\ref{thm:Th5} by taking into account the continuity of $t\mapsto u_{i}(t)\in
H$. As before we have $\Phi^{N}(-1)(t)=U_{1}^{N}(t,T)=u_{1}(t)$ and $\Phi
^{N}(1)(t)=U_{2}^{N}(t,T)=u_{2}(t)$. Moreover, the mapping $\rho
\rightarrow\Phi^{N}(\rho)(t)$ is continuous for any fixed $N\geq1$ and
$t\in\lbrack\tau,T]$ (because $U_{1}^{N}(t,\tau)=U_{2}^{N}(t,\tau)$) and
$U_{\rho\in\lbrack-1,1]}\Phi^{N}(\rho)(t)$ is connected in $H$ for each fixed
$N\geq1$ and $t\in\lbrack\tau,T]$. In particular $U_{\rho\in\lbrack-1,1])}%
\Phi^{N}(\rho)(t^{\ast})$ is connected. Indeed, we define $F:[-1,1]\rightarrow
H$ by $F(\rho):=\Phi^{N}(\rho)(t)$ and $A=[-1,1]$. If $F(A)$ were not
connected, then there would exist open sets $O_{1}$ and $O_{2}$ in $H$ with
$O_{1}\cap O_{2}=\emptyset$ such that $F(A)\cap O_{i}\neq\emptyset$, $i=1,2$,
and $F(A)\subset O_{1}\cup O_{2}$. Denote $M_{i}=\{\rho\in A:F(\rho)\in
O_{i}\}$. It is obvious that $M_{1}\cup M_{2}=A$. Moreover $M_{1}\cap
M_{2}=\emptyset$ and $M_{i}\neq\emptyset$ for $i=1,2$. Since $F$ is
continuous, it is easy to see that $V_{i}=\{\rho\in A:F(\rho)\subset O_{i}\}$
are open sets for $i=1,2$. Then $M_{i}\subset V_{i}$ and $V_{1}\cap
V_{2}=\emptyset$, which contradicts the fact that $A$ is a connected set.

Then there is $\rho_{N}\in\lbrack-1,1]$ such that $u^{N}(t^{\ast})=u_{1}%
^{N}(t,\gamma(\rho_{N}))\notin U_{1}\cup U_{2}$ and, passing to a subsequence
$\{N_{j}\}$, we prove as in Theorem {\ref{Thm:ThCON} that $u^{N_{j}}(t^{\ast
})$ converges weakly to an element $u(t^{\ast})\in K_{t^{\ast}}$ where
$u(\cdot)\in D_{\tau}(u_{\tau})$. We need to prove that the convergence is
strong. }

From
\begin{equation}
|u(t)|^{2}+\mu\int_{s}^{t}||u(r)||_{V}^{2}dr+2\alpha\int_{s}^{t}%
|u(r)|_{\beta+1}^{\beta+1}dr\leq|u(s)|^{2}+\frac{1}{\mu\lambda_{1}}\int%
_{s}^{t}|f(r)|^{2}dr\label{eq:NSDE-26BB}%
\end{equation}
and
\begin{equation}
|u^{N_{j}}(t)|^{2}+\mu\int_{s}^{t}||u^{N_{j}}(r)||^{2}dr+2\alpha\int_{s}%
^{t}\overline{F}_{N}(||u^{N_{j}}(t)||^{\beta-1})|u^{N_{j}}(r)|_{\beta
+1}^{\beta+1}dr\leq|u^{N_{j}}(s)|^{2}+\frac{1}{\mu\lambda_{1}}\int_{s}%
^{t}|f(r)|^{2}dr,\label{eq:NSDE-32BB}%
\end{equation}
for all $t>\tau$, we have
\begin{equation}
|u(t)|^{2}\leq|u(s)|^{2}+\frac{1}{\mu\lambda_{1}}\int_{s}^{t}|f(r)|^{2}%
dr\label{eq:NSDE-26BBB}%
\end{equation}
and
\begin{equation}
|u^{N_{j}}(t)|^{2}\leq|u^{N_{j}}(s)|^{2}+\frac{1}{\mu\lambda_{1}}\int_{s}%
^{t}|f(r)|^{2}dr\label{eq:NSDE-32BBB}%
\end{equation}
for all $t>s\geq\tau$. By \eqref{eq:NSDE-26BBB} and \eqref{eq:NSDE-32BBB} and
the continuity of $u^{N_{j}}(t)$ and $u(t)$ the functions $J(s)=|u(s)|^{2}%
-\frac{1}{\mu\lambda_{1}}\int_{\tau}^{s}|f(r)|^{2}dr$, $J^{N_{j}}%
(s)=|u^{N_{j}}(s)|^{2}-\frac{1}{\mu\lambda_{1}}\int_{\tau}^{s}|f(r)|^{2}dr$
are continuous and non-increasing on $[\tau,T]$.

We state that $limsupJ^{N_{j}}(t^{\ast})\leq J(t^{\ast})$. As $N_{j}%
\rightarrow\infty,$ we know that $J^{N_{j}}(t)\rightarrow J(t)$, for a.a.
$t\in(\tau,T)$. Let $t_{m}$ be a sequence such that $\tau<t_{m}\leq t^{\ast}$,
$t_{m}\rightarrow t^{\ast}$ as $m\rightarrow\infty$ and $J^{N_{j}}%
(t_{m})\rightarrow J(t_{m}),$as $m\rightarrow\infty$, for any fixed $m$.
Hence, using the continuity of J and the monotonicity of $J^{N_{j}},J$ we know
that for any $\epsilon>0$ there exist $m(\epsilon)$ and $K(\epsilon,m)$ such
that
\[
J^{N_{j}}(t^{\ast})-J(t)=J^{N_{j}}(t^{\ast})-J^{N_{j}}(t_{m})+J^{N_{j}}%
(t_{m})-J(t_{m})+J(t_{m})-J(t^{\ast})\leq|J^{N_{j}}(t_{m})-J(t_{m}%
)|+|J(t_{m})-J(t^{\ast})|\leq2\epsilon,
\]
if $N_{j}\geq K$. Thus,%
\[
limsupJ^{N_{j}}(t^{\ast})=limsup|u^{N_{j}}(t^{\ast})|^{2}-\frac{1}{\mu
\lambda_{1}}\int_{\tau}^{t^{\ast}}|f|^{2}dr\leq|u(t^{\ast})|^{2}-\frac{1}%
{\mu\lambda_{1}}\int_{\tau}^{t^{\ast}}|f|^{2}dr.
\]
Therefore, $limsup|u^{N_{j}}(t^{\ast})|\leq|u(t^{\ast})|$. Since $u^{N_{j}%
}(t^{\ast})\rightarrow u(t^{\ast})$ weakly in $H$, we have $liminf|u^{N_{j}%
}(t^{\ast})|\geq|u(t^{\ast})|$, so $u^{N_{j}}(t^{\ast})\rightarrow u(t^{\ast
})$ strongly in $H$. This gives a contradiction.
\end{proof}

\section{\textbf{Strong connectedness of the attractor for $\beta=3, $
}$\mathbf{4}$\textbf{$\alpha\mu<1$.}}

We now consider the autonomous case, i.e., $f\in H$ does not depend on $t$,
and we set $\tau=0$. Our aim now is to prove that the global compact attractor
of \eqref{eq:NSDE-1-intro}, which existence was proved in \cite[Theorem
3.1]{bib32}, is connected with respect to the strong topology of $H$ when
$\beta=3$ and $4\alpha\mu<1$.

We recall first the main results about global attractors from \cite{bib32}.

For $\beta\geq3$ we define the set
\[
\mathcal{D}(u_{0})=\{u:u\text{ is a weak global solution to
(\ref{eq:NSDE-1-intro}) with }u\left(  0\right)  =u_{0}\},
\]
which coincides with the set $\mathcal{D}_{0}(u_{0})$ defined above by
Corollary \ref{EnergyEq}.

Let $P(H)$ be the set of all non-empty sets from $H$. We define the possibly
multivalued operator $G:\mathbb{R}^{+}\times H\rightarrow P(H)$ by
\[
G(t,u_{0})=\{u(t):u\in\mathcal{D}(u_{0})\}.
\]
This operator is a strict multivalued semiflow \cite{bib32}, that is,
$G(0,u_{0})=u_{0}$ and $G(t+s,u_{0})=G(t,G(s,u_{0}))$, for all $t,s\geq
0,u_{0}\in H$.

\begin{lemma}
\label{Converg}\cite[Lemma 3.9]{bib32} Let $\beta\geq3$ and $u_{0}%
^{n}\rightarrow u_{0}$ in $H$. Then for any sequence $\{u_{n}\}_{n\geq1}$,
where $u_{n}\in\mathcal{D}(u_{0}^{n})$, there is a subsequence $\{u_{n_{k}}\}$
and $u_{0}\in\mathcal{D}(u_{0})$ such that%
\[
u_{n_{k}}\rightarrow u\text{ in }C([0,T],H)\text{ for any }T>0.
\]

\end{lemma}

Lemma \ref{Converg} implies easily that for all $t\geq0$ the operator
$G\left(  t,\text{\textperiodcentered}\right)  $ has compact values and that
for all $t\geq0$ the operator $G\left(  t,\text{\textperiodcentered}\right)  $
is upper semicontinuous, which means that for any $u_{0}\in H$ and any
neighborhood $O$ of $G(t,u_{0})$ there exists $\delta>0$ such that if
$\left\Vert v_{0}-u_{0}\right\Vert <\delta$, then $G\left(  t,v_{0}\right)
\subset O$.

The set $\mathcal{A}$ is called a global attractor for $G$ if:

\begin{enumerate}
\item $\mathcal{A}\subset G(t,\mathcal{A})$, for all $t\geq0$ (negatively invariance);

\item $dist_{H}\left(  G(t,\mathcal{A}),\mathcal{A}\right)  \rightarrow0$, as
$t\rightarrow+\infty$ (attraction property).
\end{enumerate}

If, moreover, $\mathcal{A}=G(t,\mathcal{A})$, for all $t\geq0$, it is said to
be invariant. Here, $dist_{X}$ stands for the Hausdorff semidistance defined
by $dist_{X}\left(  C,D\right)  =\sup_{c\in C}\inf_{d\in D}\left\Vert
c-d\right\Vert $, where $X$ is a Banach space.

\begin{theorem}
\cite[Theorem 3.11]{bib32} If $\beta\geq3$, then $G$ has the global compact
invariant attractor $\mathcal{A}$, which is minimal among all closed
attracting sets.
\end{theorem}

When either $\beta>3$ or $\beta=3$, $4\alpha\mu\geq1$, we have seen that the
weak solution corresponding to any $u_{0}\in H$ is unique. In such a case, the
operator $G$ becomes a single-valued semigroup and, since the space $H$ is
connected, the global attractor $\mathcal{A}$ is known to be connected
\cite{bib05}. Additionally, in \cite{bib32}, the attractor was shown to be
more regular for $\beta\geq4$. However, since the proof used $-\Delta u$ as
test function, this is unclear (see the introduction). Nevertheles, we can
give a different proof for either $\beta\in(3,5)$ or $\beta=3$, $4\alpha
\mu>1.$ For $\beta\geq5$ the question remains thus open. The proof of the
following theorem is given in the Appendix. 

\begin{theorem}
\label{AttrReg}Let either $\beta\in(3,5)$ or $\beta=3$, $4\alpha\mu>1$. Then
the global attractor $\mathcal{A}$ is bounded in $\left(  H^{2}(\Omega
)\right)  ^{3}$, which implies that it is compact in $V$ and $\left(
L^{\beta+1}(\Omega)\right)  ^{3}$. Moreover, for any bounded set $B$ from $H$
we have that%
\[
dist_{V}(G(t,B),\mathcal{A})\rightarrow0\text{ as }t\rightarrow+\infty\text{,
}%
\]
so the attraction property holds true in the stronger phase space $V.$
\end{theorem}

We recall that $V\subset\left(  L^{\beta+1}(\Omega)\right)  ^{3}$ if
$\beta\leq5$. It is known \cite{bib27} that if either $\beta\in(3,5)$ or
$\beta=3$, $4\alpha\mu>1,\ u_{0}\in V$ and $f\in H$, then there exists a
unique strong solution to problem (\ref{eq:NSDE-1-intro}). Let $S_{V}%
:\mathbb{R}^{+}\times V\rightarrow V$ be the semigroup generated by the strong
solutions. Since $S_{V}=G$ on $\mathbb{R}^{+}\times V$, the following result
follows directly from Theorem \ref{AttrReg}.

\begin{theorem}
Let either $\beta\in(3,5)$ or $\beta=3$, $4\alpha\mu>1$. Then $\mathcal{A}$ is
also a global attractor for $S_{V}.$
\end{theorem}

The last theorem gives an alternative proof of the existence of the strong
attractor from the ones in \cite{bib15}, \cite{bib16}, which are unclear. For
$\beta\geq5$ the problem remains open. Although in \cite{bib32} a conditional
result was stated for $\beta>5$, again the proof is unclear.

It remains to prove the connectedness of $\mathcal{A}$ when $\beta=3$,
$4\alpha\mu<1.$

The map $u:[0,\infty)\rightarrow H$ is said to be a trajectory of $G$ if
$u(t+s)\in G(t+s,u(s))$ for all $s,t\geq0$. Denote by $\mathcal{P}(u_{0})$ the
set of all trajectories starting at $u_{0}$. We recall that a semiflow $G $ is
called time-continuous if it is the union of continuous trajectories, that is,%
\begin{equation}
G(t,u_{0})=\{u(t):u\in\mathcal{P}(u_{0}),\ u\in C([0,\infty),H)\}.
\label{TimeCont}%
\end{equation}
Let us prove that this equality is true. In view of Lemmas \ref{Concatenation}%
, \ref{Converg} and the existence of weak solutions, the following properties
are true:

\begin{itemize}
\item[$\left(  K1\right)  $] $\mathcal{D}(u_{0})$ is non-empty for any
$u_{0}\in H.$

\item[$\left(  K2\right)  $] If $u\in\mathcal{D}(u_{0})$, then $u_{s}\left(
\text{\textperiodcentered}\right)  =u\left(  \text{\textperiodcentered
}+s\right)  \in\mathcal{D}(u_{0})$ for any $s\geq0.$

\item[$\left(  K3\right)  $] If $u,v\in\mathcal{D}(u_{0})$ are such that
$v\left(  0\right)  =u\left(  s\right)  $, where $s>0$, then%
\[
w\left(  t\right)  =\left\{
\begin{array}
[c]{c}%
u\left(  t\right)  \text{ if }0\leq t\leq s,\\
v\left(  t-s\right)  \text{ if }t\geq s,
\end{array}
\right.
\]
belongs to $\mathcal{D}(u_{0}).$

\item[$\left(  K4\right)  $] If $\{u_{n}\}$ is a sequence such that $u_{n}%
\in\mathcal{D}(u_{0}^{n})$ and $u_{n}\rightarrow u_{0}$, then there is a
subsequence $\{u_{n_{k}}\}$ and $u\in\mathcal{D}(u_{0})$ such that $u_{n_{k}%
}\left(  t\right)  \rightarrow u\left(  t\right)  $ $\forall t\geq0.$
\end{itemize}

By properties $\left(  K1\right)  -\left(  K4\right)  $ and Lemma 8 in
\cite{bib06} we get
\[
\mathcal{D}(u_{0})=\{u:u\in\mathcal{P}(u_{0}),\ u\in C([0,\infty),H)\},
\]
which implies that (\ref{TimeCont}) holds true. We observe that although Lemma
8 in \cite{bib06} is stated for complete trajectories, the proof works for
trajectories as well. Hence, our semiflow $G$ is time-continuous.

\begin{theorem}
\label{thm:Th7} Let $\beta=3$, $4\alpha\mu<1.$ Then the global attractor
$\mathcal{A}$ is connected with respect to the strong topology in $H.$
\end{theorem}

\begin{proof}
We know that the semiflow $G$ is strict and time-continuous. By Theorem
\ref{Thm:ThCON}\ we have that $G\left(  t,u_{0}\right)  $ is a connected set
in $H$ for all $t\geq0,\ u_{0}\in H$. Also, for all $t\geq0$ the operator
$G\left(  t,\text{\textperiodcentered}\right)  $ has compact values and is
upper semicontinuous. Since the space $H$ is connected, it follow from Theorem
5 in \cite{bib13} that $\mathcal{A}$ is connected with respect to the strong
topology in $H$.
\end{proof}

\bigskip

Connectedness of the attractor is an important topological property. Notably,
the Kneser property is key to applying Theorem 5 from \cite{bib13} to obtain
the connectedness of the global attractor.

\section{Appendix: proof of Theorem \ref{AttrReg}}

We will prove that the global attractor is more regular if $\beta\in(3,5)$ or
$\beta=3$ and $\mu a>1/4$. More precisely, let us check that $\mathcal{A}$ is
in fact bounded in the space $\left(  H^{2}\left(  \Omega\right)  \right)
^{3}$, and then compact in $V$ and $\left(  L^{\beta+1}\left(  \Omega\right)
\right)  ^{3}.$

\begin{lemma}
\label{EstWeakSol2}Let $\beta\in(3,5)$ or $\beta=3$ and $4\mu a>1$. Then any
weak solution to (\ref{eq:NSDE-1-intro}) with initial data such that
$\left\vert u_{0}\right\vert \leq R$ satisfies the estimate%
\begin{equation}
\left\vert u_{t}(\overline{t}+r)\right\vert ^{2}+\left\Vert u\left(
\overline{t}+r\right)  \right\Vert ^{2}+\left\vert u\left(  \overline
{t}+r\right)  \right\vert _{\beta+1}^{\beta+1}\leq D\left(  R,r\right)  ,
\label{EstNorms}%
\end{equation}
for any $r>0$ and $\overline{t}\geq0$, where $D\left(  R,r\right)  $ is such
that $D\left(  R,r\right)  \rightarrow\infty$ if $r\rightarrow0^{+}$ or
$R\rightarrow+\infty$.
\end{lemma}

\begin{proof}
Since we have uniqueness of the Cauchy problem, the following formal
calculations can be justified via Galerkin Approximations.

We prove first the result for $\overline{t}=0$. Multiplying the equation by
$Au$ we obtain that%
\[
\frac{1}{2}\frac{d}{dt}\left\Vert u\right\Vert ^{2}+\mu\left\vert
Au\right\vert ^{2}=-b(u,u,Au)-\alpha\left(  \left\vert u\right\vert ^{\beta
-1}u,Au\right)  +\left(  f,Au\right)  .
\]
By%
\[
\left\vert b(u,u,Au)\right\vert \leq C_{1}\left\Vert u\right\Vert ^{\frac
{3}{2}}\left\vert Au\right\vert ^{\frac{3}{2}}\leq\frac{\mu}{8}\left\vert
Au\right\vert ^{2}+C_{2}\left\Vert u\right\Vert ^{6},
\]%
\[
\left\vert \left(  f,Au\right)  \right\vert \leq\frac{\mu}{4}\left\vert
Au\right\vert ^{2}+\frac{1}{\mu}\left\vert f\right\vert ^{2},
\]%
\[
\alpha\left\vert \left(  \left\vert u\right\vert ^{\beta-1}u,Au\right)
\right\vert \leq\frac{\mu}{8}\left\vert Au\right\vert ^{2}+C_{3}\left\vert
u\right\vert _{2\beta}^{2\beta},
\]
where we have used inequality (9.27) in \cite{bib14}, we get%
\[
\frac{1}{2}\frac{d}{dt}\left\Vert u\right\Vert ^{2}+\frac{\mu}{2}\left\vert
Au\right\vert ^{2}\leq\frac{1}{\mu}\left\vert f\right\vert ^{2}+C_{2}%
\left\Vert u\right\Vert ^{6}+C_{3}\left\vert u\right\vert _{2\beta}^{2\beta}.
\]
Also, $\left\vert u\right\vert _{\infty}^{2}\leq C_{4}\left\Vert u\right\Vert
\left\vert Au\right\vert $ \cite{bib37} gives%
\begin{align*}
\left\vert u\right\vert _{2\beta}^{2\beta}  &  =\int_{\Omega}\left\vert
u\right\vert ^{2\beta-6}\left\vert u\right\vert ^{6}dx\leq\left\vert
u\right\vert _{\infty}^{2\beta-6}\left\vert u\right\vert _{6}^{6}\\
&  \leq C_{4}\left\Vert u\right\Vert ^{\beta+3}\left\vert Au\right\vert
^{\beta-3}\leq\frac{\mu}{4C_{3}}\left\vert Au\right\vert ^{2}+C_{5}\left\Vert
u\right\Vert ^{\frac{2(\beta+3)}{5-\beta}}.
\end{align*}
Thus,%
\begin{equation}
\frac{d}{dt}\left\Vert u\right\Vert ^{2}+\frac{\mu}{2}\left\vert Au\right\vert
^{2}\leq C_{6}(1+\left\Vert u\right\Vert ^{2q}), \label{IneqQ}%
\end{equation}
where $q=(\beta+3)/(5-\beta)$. The function $y(t)=1+\left\Vert u(t)\right\Vert
^{2}$ satisfies%
\[
y^{\prime}\leq C_{6}(1+(y(t)-1)^{q})\leq C_{6}y^{q},
\]
where we have used that $1+(y-1)^{q}\leq y^{q}$ for any $y\geq1.$ Hence,%
\[
\left(  y\left(  t\right)  \right)  ^{1-q}\geq((y(t_{0}))^{1-q}+\left(
1-q\right)  C_{6}(t-t_{0}),
\]%
\[
y(t)\leq\frac{y(t_{0})}{\left(  1+(y(t_{0}))^{q-1}\left(  1-q\right)
C_{6}(t-t_{0})\right)  ^{1/(q-1)}},
\]
which is finite if%
\[
t<t_{0}+\frac{1}{(y(t_{0}))^{q-1}(q-1)C_{6}}.
\]
Take $T^{\prime}=t_{0}+1/(2(y(t_{0}))^{q-1}(q-1)C_{6})$. Then%
\begin{equation}
\left\Vert u(t)\right\Vert ^{2}\leq2^{1/(q-1)}\left(  1+\left\Vert
u(t_{0})\right\Vert ^{2}\right)  \text{ for }t\in\lbrack t_{0},T^{\prime}].
\label{IneqnormH1}%
\end{equation}
By Lemma 2.3 in \cite{bib32} we have%
\[
\mu\int_{0}^{r}\left\Vert u(s)\right\Vert ^{2}ds\leq\left\vert u_{0}%
\right\vert ^{2}+\frac{1+r\mu\lambda_{1}}{\mu^{2}\lambda_{1}^{2}}\left\vert
f\right\vert ^{2},
\]
for $r>0$, which implies the existence of $t_{0}\in\left(  0,r\right)  $ such
that%
\[
\left\Vert u(t_{0})\right\Vert ^{2}\leq\frac{1}{\mu r}\left(  \left\vert
u_{0}\right\vert ^{2}+\frac{1+r\mu\lambda_{1}}{\mu^{2}\lambda_{1}^{2}%
}\left\vert f\right\vert ^{2}\right)  \leq D_{1}(R,r).
\]
Hence, (\ref{IneqnormH1}) implies%
\[
\left\Vert u(t)\right\Vert ^{2}\leq2^{1/(q-1)}\left(  1+\frac{1}{\mu r}\left(
\left\vert u_{0}\right\vert ^{2}+\frac{1+r\mu\lambda_{1}}{\mu^{2}\lambda
_{1}^{2}}\left\vert f\right\vert ^{2}\right)  \right)  \leq D_{2}(R,r)\text{
}\forall\ t_{0}\leq t\leq T^{\prime}\text{.}%
\]
Let $\gamma=\min\{T^{\prime},r\}$. From (\ref{IneqQ}) we have%
\begin{equation}
\sup_{t_{0}\leq t\leq\gamma}\left\Vert u(t)\right\Vert ^{2}+\frac{\mu}{2}%
\int_{t_{0}}^{\gamma}\left\vert Au\right\vert ^{2}dt\leq D_{1}(R,r)+C_{6}%
(1+\left(  D_{2}(R,r)\right)  ^{2q})r=D_{3}(R,r). \label{IneqAu}%
\end{equation}

Multiplying the equation by $u_{t}$ and using again inequality (9.27) in
\cite{bib14} we obtain%
\begin{align}
&  \left\vert u_{t}\right\vert ^{2}+\frac{\mu}{2}\frac{d}{dt}\left\Vert
u\right\Vert ^{2}+\frac{\alpha}{\beta+1}\frac{d}{dt}\left\vert u\right\vert
_{\beta+1}^{\beta+1}\nonumber\\
&  \leq-b(u,u,u_{t})+(f,u_{t})\leq C_{7}\left\Vert u\right\Vert ^{\frac{3}{2}%
}\left\vert Au\right\vert ^{\frac{1}{2}}\left\vert u_{t}\right\vert
+\left\vert f\right\vert ^{2}+\frac{1}{4}\left\vert u_{t}\right\vert
^{2}\nonumber\\
&  \leq\frac{1}{2}\left\vert u_{t}\right\vert ^{2}+C_{8}(\left\vert
f\right\vert ^{2}+\left\vert Au\right\vert ^{2}+\left\Vert u\right\Vert ^{6}).
\label{Inequt}%
\end{align}
Integrating over $\left(  t_{0},\gamma\right)  $ and using (\ref{IneqAu}) and
the embedding $V\subset\left(  L^{\beta+1}\left(  \Omega\right)  \right)  ^{3}
$ we have%
\begin{align*}
\int_{t_{0}}^{\gamma}\left\vert u_{t}\right\vert ^{2}dt  &  \leq\mu\left\Vert
u(t_{0})\right\Vert ^{2}+\frac{2\alpha}{\beta+1}\left\vert u(t_{0})\right\vert
_{\beta+1}^{\beta+1}+2C_{8}\left(  \left\vert f\right\vert ^{2}r+\frac{2}{\mu
}D_{3}(R,r)+\left(  D_{3}(R,r)\right)  ^{3}r\right) \\
&  \leq\mu D_{1}(R,r)+C_{9}\left(  D_{1}(R,r)\right)  ^{\frac{\beta+1}{2}%
}+2C_{8}\left(  \left\vert f\right\vert ^{2}r+\frac{2}{\mu}D_{3}(R,r)+\left(
D_{3}(R,r)\right)  ^{3}r\right)  =D_{4}(R,r).
\end{align*}
Hence, there is $t_{1}\in\left(  t_{0},\gamma\right)  $ such that%
\begin{equation}
\left\vert u_{t}(t_{1})\right\vert ^{2}\leq\frac{D_{4}(R,r)}{\gamma-t_{0}}%
\leq\frac{D_{4}(R,r)}{T^{\prime}-t_{0}}\leq D_{4}(R,r)2((1+D_{1}%
(R,r))^{q-1}(q-1)C_{6})=D_{5}(R,r). \label{Estut}%
\end{equation}

Further, we differentiate the equation with respect to $t$ and multiply by
$u_{t}:$%
\begin{align*}
&  \frac{1}{2}\frac{d}{dt}\left\vert u_{t}\right\vert ^{2}+\mu\left\Vert
u_{t}\right\Vert ^{2}+\alpha(\left(  \left\vert u\right\vert ^{\beta
-1}u\right)  _{t},u_{t})\\
&  =-b(u_{t},u,u_{t})-b(u,u_{t},u_{t})=b(u_{t},u_{t},u)\leq\varepsilon_{1}%
\mu\left\Vert u_{t}\right\Vert ^{2}+\frac{1}{4\varepsilon_{1}\mu}\left\vert
uu_{t}\right\vert ^{2},
\end{align*}
for $\varepsilon_{1}>0.$ Using%
\[
(\left(  \left\vert u\right\vert ^{\beta-1}u\right)  _{t},u_{t})=(\left\vert
u\right\vert ^{\beta-1}u_{t},u_{t})+\left(  \beta-1\right)  \int_{\Omega
}\left\vert u\right\vert ^{\beta-1}\left\vert u_{t}\right\vert ^{2}%
dx\geq(\left\vert u\right\vert ^{\beta-1}u_{t},u_{t})
\]
we obtain%
\begin{equation}
\frac{1}{2}\frac{d}{dt}\left\vert u_{t}\right\vert ^{2}+\mu(1-\varepsilon
_{1})\left\Vert u_{t}\right\Vert ^{2}+\alpha\left\vert \left\vert u\right\vert
^{\frac{\beta-1}{2}}\left\vert u_{t}\right\vert \right\vert ^{2}\leq\frac
{1}{4\varepsilon_{1}\mu}\left\vert uu_{t}\right\vert ^{2}. \label{Inequt2}%
\end{equation}
From (\ref{Inequt}) and%
\[
-b(u,u,u_{t})=b(u,u_{t},u)\leq\varepsilon_{2}\mu\left\Vert u_{t}\right\Vert
^{2}+\frac{1}{4\mu\varepsilon_{2}}\left\vert u\right\vert _{4}^{4},
\]
for $\varepsilon_{2}>0$, we have%
\begin{equation}
\frac{1}{2}\left\vert u_{t}\right\vert ^{2}+\frac{\mu}{2}\frac{d}%
{dt}\left\Vert u\right\Vert ^{2}+\frac{\alpha}{\beta+1}\frac{d}{dt}\left\vert
u\right\vert _{\beta+1}^{\beta+1}\leq\varepsilon_{2}\mu\left\Vert
u_{t}\right\Vert ^{2}+\frac{1}{4\mu\varepsilon_{2}}\left\vert u\right\vert
_{4}^{4}+\frac{1}{2}\left\vert f\right\vert ^{2}. \label{Inequt3}%
\end{equation}
Summing (\ref{Inequt2}) and (\ref{Inequt3}) we get%
\begin{align*}
&  \frac{d}{dt}\left(  \frac{1}{2}\left\vert u_{t}\right\vert ^{2}+\frac{\mu
}{2}\left\Vert u\right\Vert ^{2}+\frac{\alpha}{\beta+1}\left\vert u\right\vert
_{\beta+1}^{\beta+1}\right)  +\frac{1}{2}\left\vert u_{t}\right\vert ^{2}%
+\mu(1-\varepsilon_{1}-\varepsilon_{2})\left\Vert u_{t}\right\Vert ^{2}%
+\alpha\left\vert \left\vert u\right\vert ^{\frac{\beta-1}{2}}\left\vert
u_{t}\right\vert \right\vert ^{2}\\
&  \leq\frac{1}{4\varepsilon_{1}\mu}\left\vert uu_{t}\right\vert ^{2}+\frac
{1}{4\mu\varepsilon_{2}}\left\vert u\right\vert _{4}^{4}+\frac{1}{2}\left\vert
f\right\vert ^{2}\leq\frac{1}{4\varepsilon_{1}\mu}\left\vert uu_{t}\right\vert
^{2}+K_{1}(\varepsilon_{2})\left(  1+\left\vert u\right\vert _{\beta+1}%
^{\beta+1}\right)  .
\end{align*}
If $\beta=3$, the condition $\alpha\mu4>1$ implies that
\[
\alpha\left\vert \left\vert u\right\vert ^{\frac{\beta-1}{2}}\left\vert
u_{t}\right\vert \right\vert ^{2}-\frac{1}{4\varepsilon_{1}\mu}\left\vert
uu_{t}\right\vert ^{2}=\left(  \alpha-\frac{1}{4\varepsilon_{1}\mu}\right)
\left\vert uu_{t}\right\vert ^{2}\geq0
\]
for $\varepsilon_{1}$ close enough to $1$. If $\beta>3$, then by Young's
inequality there is $K_{2}(\varepsilon_{1})>0$ such that%
\[
\alpha\left\vert \left\vert u\right\vert ^{\frac{\beta-1}{2}}\left\vert
u_{t}\right\vert \right\vert ^{2}-\frac{1}{4\varepsilon_{1}\mu}\left\vert
uu_{t}\right\vert ^{2}\geq\frac{\alpha}{2}\left\vert \left\vert u\right\vert
^{\frac{\beta-1}{2}}\left\vert u_{t}\right\vert \right\vert ^{2}%
-K_{2}(\varepsilon_{1})\left\vert u_{t}\right\vert ^{2}.
\]
Hence,%
\[
\frac{d}{dt}\left(  \frac{1}{2}\left\vert u_{t}\right\vert ^{2}+\frac{\mu}%
{2}\left\Vert u\right\Vert ^{2}+\frac{\alpha}{\beta+1}\left\vert u\right\vert
_{\beta+1}^{\beta+1}\right)  \leq K_{1}(\varepsilon_{2})\left(  1+\left\vert
u\right\vert _{\beta+1}^{\beta+1}\right)  +K_{2}(\varepsilon_{1})\left\vert
u_{t}\right\vert ^{2},\ \forall\ t\geq t_{1},
\]
if we choose $\varepsilon_{1}$ close enough to $1$ and $\varepsilon_{2}$ small
enough. Let $y(t)=\frac{1}{2}\left\vert u_{t}\right\vert ^{2}+\frac{\mu}%
{2}\left\Vert u\right\Vert ^{2}+\frac{\alpha}{\beta+1}\left\vert u\right\vert
_{\beta+1}^{\beta+1}$. Then using Gronwall's lemma, the embedding
$V\subset\left(  L^{\beta+1}(\Omega)\right)  ^{3}$ and (\ref{IneqAu}),
(\ref{Estut}) we obtain for some constants $K=K(\varepsilon_{1},\varepsilon
_{2}),\ D_{6}(R,r)$ that%
\[
y(t)\leq\left(  y(t_{1})+1\right)  e^{K(t-t_{1})}\leq D_{6}(R,r)e^{K(t-t_{1}%
)}\ \forall t\geq t_{1}.
\]
Thus, in particular,%
\[
y(r)\leq D_{6}(R,r)e^{K(r-t_{1})}=D_{7}(R,r),
\]
which proves (\ref{EstNorms}) for $\overline{t}=0.$

For an arbitrary $\overline{t}\geq0$ by Lemma 2.3 in \cite{bib32} we make use
of the estimate%
\[
\left\vert u(t)\right\vert ^{2}\leq\left\vert u(0)\right\vert ^{2}%
+\frac{\left\vert f\right\vert ^{2}}{\mu^{2}\lambda_{1}^{2}}\leq R^{2}%
+\frac{\left\vert f\right\vert ^{2}}{\mu^{2}\lambda_{1}^{2}}=R_{1}%
^{2},\ \forall t\geq0.
\]
Defining $v(t)=u(t+\overline{t})$, then
\[
\frac{1}{2}\left\vert v_{t}\left(  r\right)  \right\vert ^{2}+\frac{\mu}%
{2}\left\Vert v(r)\right\Vert ^{2}+\frac{\alpha}{\beta+1}\left\vert v\left(
r\right)  \right\vert _{\beta+1}^{\beta+1}\leq D_{7}\left(  R_{1},r\right)
=D_{8}(R,r),
\]
which gives (\ref{EstNorms}).
\end{proof}

\bigskip

As a consequence of this lemma and the compact embedding $V\subset H$ we
obtain the following result.

\begin{corollary}
\label{CompactOperator}Let $\beta\in(3,5)$ or $\beta=3$ and $4\mu a>1$. Then
for any $r>0$ the map $u_{0}\mapsto S\left(  r,u_{0}\right)  $ maps bounded
subsets of $H$ onto bounded subsets of $V\cap L^{\beta+1}\left(
\Omega\right)  $. Hence, $S\left(  r\right)  $ is a compact operator, i.e., it
maps bounded subsets of $H$ onto relatively compact ones.
\end{corollary}

\begin{lemma}
\label{EstAu}Let $\beta\in(3,5)$ or $\beta=3$ and $4\mu a>1$. Then any weak
solution to (\ref{eq:NSDE-1-intro}) with initial data such that $\left\vert
u_{0}\right\vert \leq R$ satisfies the estimate%
\[
\left\vert Au\left(  r\right)  \right\vert \leq K\left(  R,r\right)  ,
\]
for any $r>0$, where $K\left(  R,r\right)  $ is such that $K\left(
R,r\right)  \rightarrow\infty$ if $r\rightarrow0^{+}$ or $R\rightarrow+\infty$.
\end{lemma}

\begin{proof}
By Proposition 9.2 in \cite{bib14} we have%
\[
\left\vert B(u,u)\right\vert \leq d_{1}\left\Vert u\right\Vert ^{\frac{3}{2}%
}\left\vert Au\right\vert ^{\frac{1}{2}}\leq\frac{\mu}{4}\left\vert
Au\right\vert +d_{2}\left\Vert u\right\Vert ^{3}.
\]
Using the Gagliardo-Nirenberg inequality and $\beta<5$ we find that%
\begin{align*}
\alpha\left\vert \left\vert u\right\vert ^{\beta-1}u\right\vert  &
=\alpha\left\vert u\right\vert _{2\beta}^{\beta}\leq d_{3}\left\vert
Au\right\vert ^{\frac{3\left(  \beta-1\right)  }{\beta+7}}\left\vert
u\right\vert _{\beta+1}^{\frac{\beta^{2}+4\beta+3}{\beta+7}}\\
&  \leq\frac{\mu}{4}\left\vert Au\right\vert +d_{4}\left\vert u\right\vert
_{\beta+1}^{\frac{\beta^{2}+4\beta+3}{10-2\beta}}.
\end{align*}

Hence,
\[
\frac{\mu}{2}\left\vert Au\left(  r\right)  \right\vert \leq\left\vert
u_{t}(r)\right\vert +d_{2}\left\Vert u(r)\right\Vert ^{3}+d_{4}\left\vert
u(r)\right\vert _{\beta+1}^{\frac{\beta^{2}+4\beta+3}{10-2\beta}}+\left\vert
f\right\vert ,
\]
so the result follows by applying Lemma \ref{EstWeakSol2}.
\end{proof}

\bigskip

We are now in position of proving the regularity of the global attractor.

\begin{theorem}
Let $\beta\in(3,5)$ or $\beta=3$ and $4\mu\alpha>1$. Then the global attractor
$\mathcal{A}$ is bounded in $\left(  H^{2}(\Omega)\right)  ^{3}$, and then
compact in $V$ and $\left(  L^{\beta+1}(\Omega)\right)  ^{3}$. Moreover,
\begin{equation}
dist_{V}(S(t,B),\mathcal{A})\rightarrow0\text{ as }t\rightarrow+\infty,
\label{ConvergVLbeta}%
\end{equation}
for any $B$ bounded in $H$.
\end{theorem}

\begin{proof}
Since the global attractor is invariant, $\mathcal{A}=S(r,\mathcal{A})$ for
$r>0$, so $\mathcal{A}$ is bounded in $\left(  H^{2}(\Omega)\right)  ^{3}$ by
Lemma \ref{EstAu}. The compact embeddings $H^{2}(\Omega)\subset H^{1}(\Omega
)$, $H^{2}(\Omega)\subset L^{\beta+1}(\Omega)$ imply the compactness of the
attractor in $V$ and $\left(  L^{\beta+1}(\Omega)\right)  ^{3}$.

Let $B_{0}$ be an absorbing ball. By Lemma \ref{EstAu} the set $B_{1}%
=S(r,B_{0})$ is bounded in $\left(  H^{2}(\Omega)\right)  ^{3}$ and%
\[
S(t,B)=S(r,S(t-r,B))\subset B_{1}\text{ for }t\geq t_{0}(B)\text{.}%
\]
From here it is easy to deduce (\ref{ConvergVLbeta}).
\end{proof}

\bigskip

\textbf{Acknowledgments}

The second atuthor has been supported by Ministerio de Ciencia, Innovaci\'{o}n
y Universidades, AEI and FEDER, project PID-156228NB-I00

\end{document}